\PassOptionsToPackage{table}{xcolor}
\documentclass[hidelinks,onefignum,onetabnum]{siamart250211}

\usepackage{lipsum}
\usepackage{amsfonts}
\usepackage{graphicx}
\usepackage{epstopdf}
\usepackage{algorithmic}
\ifpdf
  \DeclareGraphicsExtensions{.eps,.pdf,.png,.jpg}
\else
  \DeclareGraphicsExtensions{.eps}
\fi

\newsiamremark{remark}{Remark}
\newsiamremark{hypothesis}{Hypothesis}
\crefname{hypothesis}{Hypothesis}{Hypotheses}
\newsiamthm{claim}{Claim}
\newsiamremark{fact}{Fact}
\crefname{fact}{Fact}{Facts}

\headers{A flow-rate-conserving CNN-based DDM for blood flow simulations}{A flow-rate-conserving CNN-based DDM for blood flow simulations}

\title{A flow-rate-conserving CNN-based domain decomposition method for blood flow simulations\thanks{Submitted to the editors \today.
\funding{This work was performed as part of the Helmholtz School for Data Science in
Life, Earth, and Energy (HDS-LEE) and received funding from the Helmholtz Association of German
Research Centres.}}}

\author{
  Simon Klaes\footnotemark[2] \and
  Axel Klawonn\footnotemark[2] \footnotemark[3] \and
  Natalie Kubicki\footnotemark[2] \footnotemark[4] \and
  Martin Lanser\footnotemark[2] \footnotemark[3] \and
  Kengo Nakajima\footnotemark[5] \and
  Takashi Shimokawabe\footnotemark[5] \and
  Janine Weber-Hamacher\footnotemark[2] \footnotemark[3]
}

\usepackage{amsopn}

\ifpdf
\hypersetup{
  pdftitle={
  },
  pdfauthor={Simon Klaes, Axel Klawonn, Natalie Kubicki,  Martin Lanser, Kengo Nakajima, Takashi Shimokawabe, Janine Weber-Hamacher  }
}
\fi

\usepackage{ifthen}
\newcommand{\mycomment}[1]{
\ifthenelse{\isodd{\value{page}}}{
\normalmarginpar%
\marginpar{\tiny {#1}}
}{
\reversemarginpar%
\marginpar{\tiny {#1}}
}}

\def\KL{\color{black}}  

\usepackage{comment}

\usepackage{tikz}
\usetikzlibrary{positioning}
\usetikzlibrary{calc} 
\definecolor{myblue}{HTML}{6595FC}
\definecolor{royalblue}{HTML}{0C3DC2}

\usepackage{pgfplots}

\usepackage{amsmath, amssymb}
\usepackage{bm}

\usepackage{booktabs}
\usepackage{multirow}
\usepackage[table]{xcolor}
\usepackage{diagbox}

\usepackage{CJKutf8}

\usepackage{ulem}
\usepackage{ifthen}

\crefname{figure}{Fig.}{Figs.}

\begin{document}

\maketitle

\renewcommand{\thefootnote}{\fnsymbol{footnote}}

\footnotetext[2]{Department of Mathematics and Computer Science, University of Cologne, Germany (\email{s.klaes@uni-koeln.de}, \email{axel.klawonn@uni-koeln.de}, \email{nkubicki@uni-koeln.de}, 
  \email{martin.lanser@uni-koeln.de}, and \email{janine.weber@uni-koeln.de}), url: \url{http://www.numerik.uni-koeln.de}} 
\footnotetext[3]{Center for Data and Simulation Science, University of Cologne, Germany, url: \url{http://www.cds.uni-koeln.de}}
\footnotetext[4]{Institute for Advanced Simulation  (IAS-2), Forschungszentrum Jülich, Germany}
\footnotetext[5]{Information Technology Center, University of Tokyo, Japan 
  (\email{nakajima@cc.u-tokyo.ac.jp} and \email{shimokawabe@cc.u-tokyo.ac.jp})} 

\begin{abstract} This work aims to predict blood flow with non-Newtonian viscosity in stenosed arteries using convolutional neural network (CNN) surrogate models. 
    An alternating Schwarz domain decomposition method is proposed which uses CNN-based subdomain solvers.
     A universal subdomain solver (USDS) is trained on a single, fixed geometry and then applied for each subdomain solve in the Schwarz method.
    Results for two-dimensional stenotic arteries of varying shape and length for different inflow conditions are presented and statistically evaluated. 
One key finding, when using a limited amount of training data, is that incorporating a physics-aware constraint, as, in our case,
flow rate conservation, into the USDS improves the prediction accuracy and convergence behavior of the Schwarz method compared to a purely data-driven USDS.
As the USDS is a data-driven, inexact subdomain solver, admissible parameter ranges for the geometry and inflow configurations must be defined and tested.
\end{abstract} 


\section{Introduction}
\label{sec:intro}
Computational fluid dynamics (CFD) is an important tool for the
in-silico investigation of biomedical processes \cite{CardiovascularMathematics}.
While the classical CFD approaches can provide accurate predictions, they become computationally expensive with an increasing 
mesh resolution. With the  
growing success 
of data-driven methods, deep learning techniques have become popular to learn complex nonlinear mappings directly from data generated by, for example, conventional high-fidelity solvers \cite{review_ml_cfd_wang}. While these  data-driven approaches can demonstrate impressive acceleration factors,  they also face several significant challenges \cite{McGreivy2024_weakbaseline}.
As a drawback, the trained surrogate models typically lack transferability, being trained for specific geometries or model parameters, making them perform poorly for different configurations. 
Additionally, purely data-driven models do not take into account underlying physical laws and often struggle with physical consistency, especially when extrapolating beyond the training data set.
Recent research tries to address these challenges. For example, the research field scientific machine learning (SciML)  combines data-driven and physics-aware approaches to improve the generalization capabilities of deep learning-based surrogate models \cite{USDpt:SciML}. In particular, SciML combines knowledge from classic numerical solvers, as, for example, domain decomposition methods, with experience in machine learning techniques with the aim to develop fast and robust hybrid methods tailored for a specific given application; see \cite{article_dd_ml_survey} for further details in the context of domain decomposition methods. 

\begin{figure}[ht!]
    \resizebox{\textwidth}{!}{ 
\includegraphics[scale=0.5]{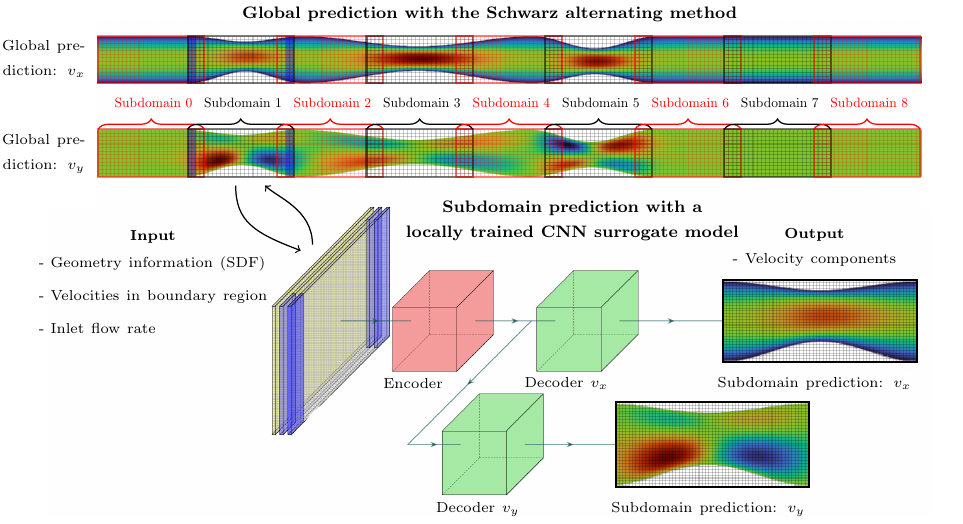}
}
\vspace{-15pt}
\caption{Visualization of the domain decomposition-based framework using a locally trained CNN for predicting stationary blood flow in stenosed arteries. 
The geometry is decomposed into over\-lapp\-ing subdomains of fixed size. 
For the global flow prediction an alternating Schwarz method is applied, where a locally trained CNN acts as an inexact subdomain solver and velocity fields are iteratively updated.  
The CNN consists of a shared encoder and separate-path decoders, taking as input geometry information and velocities from a boundary region of the subdomain to predict the velocity components within each subdomain.
To enforce flow rate conservation throughout the domain, a physics-aware constraint layer is integrated into the CNN that scales the locally predicted velocity field such that the flow rate in each cross-section corresponds to a given inlet flow rate.}
\label{fig:overview}
\end{figure}

In the present work, we focus on the development of surrogate models for solving the forward problem arising in macroscopic blood flow simulations. 
Particulary, we are interested in how variations in arterial geometries, such as stenoses, and different inflow conditions influence the resulting flow fields in arteries of arbitrary length. 
For the development of a surrogate model, we build on \cite{guo_cnn, eichinger_cnn}, utilizing CNNs as data-driven surrogate models for predicting stationary flow fields.
CNNs  are a class of deep neural networks specifically designed for processing structured grid data such as images; see \cite{lecun:1989:CNN,Goodfellow:2016:DL} for details.
 When training CNNs for CFD applications, the key idea is to 
 train a CNN to approximate a solution operator, which maps an image representation of the specified inputs onto an image representation of the outputs, 
 that is, velocity components. However, two significant limitations 
 are the necessity of large training datasets  and the limited generalizability 
of the trained models to varying domain sizes, as the input image size is fixed by the network architecture and cannot be changed directly.
  To adress these limitations, in this work, we investigate the combination of CNNs with established domain decomposition strategies \cite{book_domaindecomposition, dd_quarteroni}.
Instead of learning the flow dynamics for the entire geometry domain at once, the global problem is decomposed into smaller, rectangular subdomains each capturing a partial image of the geometry; see \cref{fig:overview}. 
The CNN is then trained to predict the flow field for the local problem within a single subdomain using as inputs geometry information decoded by a signed distance function (SDF) and velocities in a certain boundary region of the subdomain. 
Based on our specific problem assumptions, we investigate the integration of a physics-aware constraint layer into the CNN that scales the predicted velocity fields column-wise such that the flow rate in each cross-section matches a provided inlet flow rate, strictly enforcing flow rate conservation within each subdomain. 
By reducing the size of the problem through the decomposition, the CNN has to learn local and thus simpler dynamics compared to considering the entire geometry at once. Now, in order to propagate information between subdomains, a coupling in form of overlaps and exchanged boundary conditions is considered. 
Hence, an alternating Schwarz strategy is implemented in order to enable the global\footnote{For those readers familiar with two-level domain decomposition preconditioners, we note that here the word {\it global} does not relate to a coarse problem but to the original problem.} prediction to converge to an approximation of the flow field for the specific global boundary conditions. 
As the CNN thus acts as an inexact, data-driven subdomain solver, it is necessary to define admissible parameter ranges for the geometry and inflow
conditions based on the training data and corresponding tests.
The key advantage of our approach is that the CNN does not need to be retrained for longer arterial geometries.

In \cite{article_rana_cnn_dd}, a CNN-based approach was introduced that extends flow predictions to domains of arbitrary sizes using, in contrast to our work, a non-overlapping domain decomposition strategy.
A closely related work to our described procedure, based on overlapping subdomains, is introduced in \cite{Wang_2022}. 
There, the underlying network architecture for the local flow prediction uses physics-informed neural networks, in which physical constraints are integrated into the loss function \cite{RAISSI2019686}.
The authors demonstrate scalability to different geometries, allowing locally trained networks to be used for various geometry boundary conditions \cite{feeney2023breakingboundariesdistributeddomain, Wang_2022}.
In \cite{new_paper_schwarz_pde_neural_operators}, the authors present a theoretical framework for using locally trained neural operators within an additive Schwarz method to solve elliptic PDEs in arbitrary domains.
In contrast, our work applies to a specific application, focusing on blood flow 
and emphasizing the importance of preserving conservation laws in the local model and defining admissible parameter ranges to ensure reliable global solutions.
For a recent survey on combining domain decomposition methods with machine learning; see \cite{article_dd_ml_survey}. Our primary objective in this paper is to provide a detailed pipeline and a proof of concept,  discussing the strengths and weaknesses of the proposed approach, 
referred to in the following as CNN-Schwarz-Flow. We do not claim to have the most advanced network architecture for the local problem or the most efficient code for reducing the training time, but we rather focus on the methodology itself.
 
The remainder of the paper is organized as follows. In \cref{sec:methods}, we introduce the non-Newtonian flow problem modeling the blood flow in an artery. 
 In \cref{sec:new_algorithm_dd_cnn}, we describe  in detail the 
 CNN-Schwarz-Flow approach. 
First, we present the CNN architecture for the local problem (see \cref{sec:single_image_prediction}) and describe how a physics-aware hard constraint based on 
 flow rate conservation can be integrated into the model. Next, the algorithm for the global prediction (see \cref{sec:global_image_prediction}) based on an alternating Schwarz strategy is discussed. 
 In \cref{sec:setup}, the data generation process and the training procedure are described.
  The results of applying the 
  CNN-Schwarz-flow approach using a flow-rate-conserving and a purely data-driven CNN are empirically analyzed and evaluated in \cref{sec:results}. 
   Our conclusions are summarized in 
  \cref{sec:conclusion}.

\section{Fluid flow problem}
\label{sec:methods}

As our data-driven surrogate models learn from  
simulation data and can only capture physics described by the data, the underlying model assumptions used for the generation of the training data are described in detail.

\begin{figure}[ht!]
    \resizebox{\textwidth}{!}{ 
\begin{tikzpicture}
  \node[anchor=south west,inner sep=0] (image) at (0,0) 
    {\includegraphics[width=\textwidth,height=0.045\textheight]{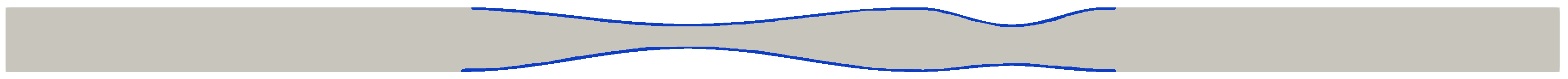}};
  
\draw[black, <->] ([xshift=0.07cm, yshift=0.1cm]image.north west) -- ([xshift=3.98cm, yshift=0.1cm]image.north west)
    node[midway, above] {$L_{\text{inlet}}$};

\draw[royalblue, line width=1pt] ([xshift=0.07cm, yshift=-0.1cm]image.north west) -- ([xshift=3.98cm, yshift=-0.1cm]image.north west);

\draw[royalblue, line width=1pt] ([xshift=0.07cm, yshift=-0.85cm]image.north west) -- ([xshift=3.98cm, yshift=-0.85cm]image.north west);

\draw[royalblue, line width=1pt] ([xshift=-0.08cm,yshift=-0.11cm]image.north east) -- ([xshift=-3.9cm, yshift=-0.11cm]image.north east);
\draw[royalblue, line width=1pt] ([xshift=-0.08cm,yshift=-0.85cm]image.north east) -- ([xshift=-3.9cm,yshift=-0.85cm]image.north east);

\draw[black, <->] ([xshift=-0.05cm, yshift=0.1cm]image.north east) -- ([xshift=-4cm, yshift=0.1cm]image.north east)
    node[midway, above]  {$L_{\text{outlet}}$};

\draw[black, <->] ([xshift=-9.02cm, yshift=0.1cm]image.north east) -- ([xshift=-4.02cm, yshift=0.1cm]image.north east)
    node[midway, above] {$L_{\text{stenotic}}$};

\draw[black, <->] ([xshift=-0.2cm, yshift=-0.08cm]image.north west) -- ([xshift=-0.2cm, yshift=-0.87cm]image.north west)
    node[midway, rotate=90, yshift=0.4cm] {$d_{\text{artery}}$};
    
\draw[red, line width=2pt] ([xshift=-0.1cm, yshift=-0.12cm]image.north east) -- ([xshift=-0.1cm, yshift=-0.84cm]image.north east);
   
\draw[violet, line width=2pt] ([xshift=0.07cm, yshift=-0.12cm]image.north west) -- ([xshift=0.07cm, yshift=-0.84cm]image.north west);

 \draw[violet, <-, line width=0.6pt]   	 ([xshift=0.2cm, yshift=-0.28cm]image.north west)  --([xshift=1.1cm, yshift=-0.5cm]image.north west) 
 node[pos=1, right] {$\Gamma_{D_{\text{inlet}}}$};

  \draw[royalblue, <-, line width=0.6pt]   	 ([xshift=3.98cm, yshift=-0.2cm]image.north west)  --([xshift=5.1cm, yshift=-0.9cm]image.north west) 
 node[pos=1, right] {$\Gamma_{D_{\text{walls}}}$};
  \draw[royalblue, <-, line width=0.6pt]   	 ([xshift=4.5cm, yshift=-0.85cm]image.north west)  --([xshift=5.1cm, yshift=-0.9cm]image.north west) ;

 \draw[red, <-, line width=0.6pt]   	([xshift=-0.17cm, yshift=-0.5cm]image.north east) -- ([xshift=-1cm, yshift=-0.52cm]image.north east)
 node[pos=1, left] {$\Gamma_N$};
  	
\node[font=\large, xshift=0.5cm] at ($(image.south west)!0.5!(image.north east)$) {$\Omega$};

\draw[->, thick, black] ([xshift=0.07cm, yshift=0.1cm]image.south west) -- ([xshift=0.4cm, yshift=0.1cm]image.south west) node[right,yshift=-2.8] {$x$};
\draw[->, thick, black] ([xshift=0.07cm, yshift=0.1cm]image.south west) -- ([xshift=0.07cm, yshift=0.45cm]image.south west) node[right,yshift=-0.5] {$y$};
    
\end{tikzpicture}
}
\vspace{-15pt} 
\caption{Sample arterial geometry with domain $\Omega \subset \mathbb{R}^{2}$. The maximum diameter $d_{\text{artery}} = 1$ [mm] and total length $L = L_{\text{inlet}} + L_{\text{stenotic}} + L_{\text{outlet}} = 2.4$ [cm] are  assumed to be constant for all geometries with $L_{\text{inlet}} = L_{\text{outlet}} = 0.7$ [cm]. The specific form of the stenotic region $L_{\text{stenotic}}=1$ [cm] can differ. For detailed information regarding the geometry generation; see \cref{geometry_fluidmodel}. At the inlet $\Gamma_{D_{\text{inlet}}}$, a parabolic velocity profile is prescribed. At the walls $\Gamma_{D_{\text{walls}}}$, no-slip conditions are set, and at the outflow $\Gamma_N$, the stress-free boundary condition is applied. The origin is located at $(0,0)$.}
\label{geometry_methods}
\end{figure}

\paragraph{\textbf{Fluid Model}}
Our fluid of interest is blood, which is in general a complex fluid exhibiting various rheological phenomena  \cite{Bessonov2016_MethodsBloodFlow, Fasano2017}. 
Here, blood flow is modeled on a macroscopic scale, treating blood as a continuum. 
Specifically, blood is assumed to be a single phase, homogeneous, incompressible, isothermal, and shear-thinning fluid under stationary flow conditions in arteries 
with  rigid vessel walls. Under these assumptions, the governing equations reduce to the mass and momentum conservation laws. In the Eulerian frame, for a domain $\Omega \subset \mathbb{R}^d$ with $d$ spatial dimensions, and with boundary $\Gamma = \partial \Omega \subset \mathbb{R}^{d-1}$ decomposed into non-overlapping Neumann ($\Gamma_N$) and Dirichlet ($\Gamma_D$) boundaries, the resulting equations are \cite{CardiovascularMathematics}
\begin{subequations}
\label{eq:gnf_pde}
\begin{align}
\bm{\nabla} \cdot \bm{v} &= 0 \quad\quad\quad\, \text{in}  \,\, \Omega, \label{eq:divergendefree} \\
\rho (\bm{v}\cdot \bm{\nabla})\bm{v} + \bm{\nabla}p - \bm{\nabla} \cdot \bm{\tau} &= \bm{f} \quad\quad\quad \text{in}  \,\, \Omega, \\
\bm{v} &= \bm{g} \quad\quad\quad\hspace{0.5pt} \text{on} \,\, \Gamma_D, \\
(-p \bm{I} + \bm{\tau}) \bm{n} &= 0 \quad\quad\quad\, \text{on} \,\, \Gamma_N.
\label{eq:bb}
\end{align}
\end{subequations}
Here, $\rho$ is the constant fluid density, $p$ and $\bm{v}$  are the unknown pressure and velocity fields, $\bm{\tau}$ is the viscous stress tensor, and $\bm{f}$ denotes external forces, which are omitted ($\bm{f} = \bm{0}$). Furthermore, $\bm{g}$ is a function specifying the Dirichlet boundary conditions, and in the stress-free Neumann boundary condition \cref{eq:bb}, $\bm{I}$ is the identity tensor, and $\bm{n}$ is the normal vector along $\Gamma_N$.
We constrain ourselves to two-dimensional  arterial  geometries, which are constructed as channel-like domains with varying diameters; see \cref{geometry_methods}. 
The Dirichlet boundary is divided as $\Gamma_D = \Gamma_{D_{\text{inlet}}} \cup \Gamma_{D_{\text{walls}}}$. At $ \Gamma_{D_{\text{inlet}}}$, a parabolic profile is prescribed,  in which the maximum velocity $v^{\max}_{\text{inlet}}$ can vary to consider different flow rates. At the walls, no-slip conditions are applied, that is,
\begin{subequations} 
\begin{align}
 \bm{v} &= \begin{pmatrix}
\frac{4 \, v^{\max}_{\text{inlet}}\, y\,(d_{\text{artery}}-y)}{d_{\text{artery}}^2} \, , \quad 0
\end{pmatrix}^{\top} \quad\hspace{0.5pt} \text{on} \,\, \Gamma_{D_{\text{inlet}}} \,\quad \text{with} \,\quad v^{\max}_{\text{inlet}} \in [0.02,0.6] \,[\frac{\text{m}}{\text{s}}], \label{parabolic} \\
\bm{v} &= \bm{0} \quad\hspace{108pt} \text{on} \,\, \Gamma_{D_{\text{walls}}}.
\label{eq:parabolic}
\end{align}
\end{subequations} As necessary closure model for the stress tensor $\bm{\tau}$, a generalized Newtonian constitutive relation with non‑Newtonian viscosity is assumed, given by
\begin{equation}
\bm{\tau} = \eta(\dot{\gamma})  \left( \bm{\nabla} \bm{v} + (\bm{\nabla} \bm{v})^{\top} \right),
\end{equation}
where the viscosity $\eta(\dot{\gamma}): \mathbb{R}^{+} \rightarrow \mathbb{R}^{+} \setminus \{0\}$ is a scalar function depending on the shear rate $\dot{\gamma}$ \cite{Sequeira2018}. 
Here, a Carreau model is assumed for the viscosity, that  is, 
\begin{align}\label{eq:viscosityModel}
    \frac{\eta(\dot{\gamma})}{\eta_{\text{ref}}} = \eta_{\infty} + (\eta_0 - \eta_{\infty}) \left( 1 + (\lambda \dot{\gamma})^2 \right)^{\frac{n-1}{2}},
\end{align}
where \(\eta_0\) and \(\eta_{\infty}\) represent the zero-shear-rate and infinite-shear-rate viscosities, \(\lambda\) is a time parameter, \(n\) is the power law index, and the model is divided by a reference viscosity \(\eta_{\text{ref}}\) \cite{irgens_rheology}. 
The values of these parameters are fitted to data provided by the Institute for Advanced Simulation  (IAS-2) at Forschungszentrum Jülich, which employed a virtual rheometer using dissipative particle dynamics simulations to obtain viscosity data for blood under
certain conditions \cite{rheometer}; see \cref{tab:viscosity_params}.
\begin{table}[htbp]
\centering
\small 
\renewcommand{\arraystretch}{1.4}  
\begin{tabular}{l c c c c c c}
\toprule[1.0pt]
\textbf{Blood parameters} & $\rho$ [kg/m$^3$] & $\eta_{\infty}$ & $\eta_0$ & $\eta_{\text{ref}}$ [Pa·s] & $n$ & $\lambda$ [s] \\
\midrule
\textbf{Values} & 1000 & 3.3707 & 230.6330 & 0.0012 & 0.45 & -300.0 \\
\bottomrule[1.0pt]
\end{tabular}
\vspace{7pt}
\caption{Specified blood model parameters. The viscosity parameters are fitted to virtual rheometer data assuming blood under specific conditions  (Hematocrit (volume percentage of red blood cells) $= 45$\%, Temperatur $= 35^{\circ}$C, healthy red blood cells).}
\label{tab:viscosity_params}
\end{table}

\begin{figure}[ht!]
    \vspace{-7pt}
    \resizebox{\textwidth}{!}{ 
\begin{tikzpicture}
    \usetikzlibrary{decorations.pathreplacing}
    
           \node[anchor=south west,inner sep=0] (image) at (0,2) 
        {\includegraphics[width=\textwidth]{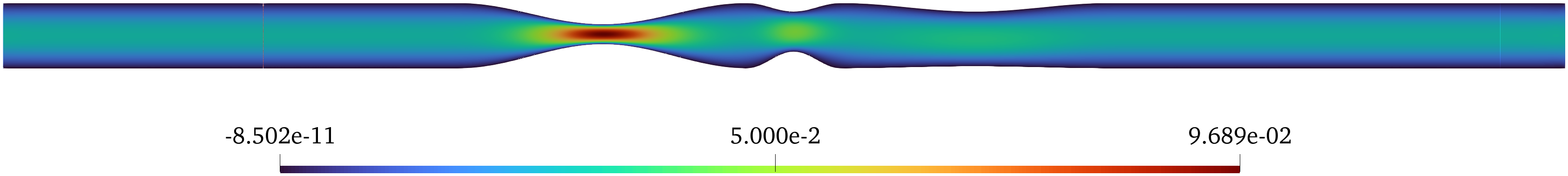}}; 
 
    \node[anchor=south west,inner sep=0] (image2) at (0,0) 
        {\includegraphics[width=\textwidth]{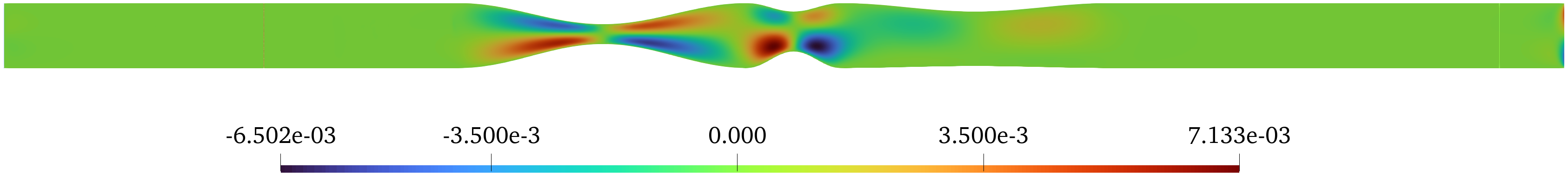}};
       
\draw[blue, line width=2pt] ([xshift=0.02cm, yshift=-0.03cm]image.north west) -- ([xshift=0.02cm, yshift=-0.58cm]image.north west)  node[pos=0, above] {$x=0$};    

\draw[green!70!black, line width=2pt] ([xshift=2.17cm, yshift=-0.03cm]image.north west) -- ([xshift=2.17cm, yshift=-0.58cm]image.north west) node[pos=0, above] {$x=x_s$};   

\draw[red, line width=2pt] ([xshift=12.48cm, yshift=-0.03cm]image.north west) -- ([xshift=12.48cm, yshift=-0.58cm]image.north west) node[pos=0, above] {$x=x_e$};

\draw[blue, line width=2pt] ([xshift=0.02cm, yshift=-0.03cm]image2.north west) -- ([xshift=0.02cm, yshift=-0.58cm]image2.north west)   node[pos=0, above] {$x=0$};    

\draw[green!70!black, line width=2pt] ([xshift=2.17cm, yshift=-0.03cm]image2.north west) -- ([xshift=2.17cm, yshift=-0.58cm]image2.north west)  node[pos=0, above] {$x=x_s$};   

\draw[red, line width=2pt] ([xshift=12.48cm, yshift=-0.03cm]image2.north west) -- ([xshift=12.48cm, yshift=-0.58cm]image2.north west) node[pos=0, above] {$x=x_e$};  
       
    \node[above] at (6.4, 0.4) {\scalebox{.5}{\large Velocity $y$-component $v_y$ $[\frac{\text{m}}{\text{s}}]$}}; 
    \node[above] at (6.4, 2.4) {\scalebox{.5}{\large Velocity $x$-component $v_x$ $[\frac{\text{m}}{\text{s}}]$}};

\end{tikzpicture}
}
\vspace{-12pt}
\caption{Visualization of computed velocity fields using $v^{\max}_{\text{inlet}}=0.03 \,[\frac{\text{m}}{\text{s}}]$ for a sample arterial geometry. The prescribed parabolic profile at the inlet ($x=0$) flattens slightly in accordance with the assumed generalized Newtonian model, reaching a fully developed state in a sufficiently long, straight channel segment.  In the cross-sections at $x_s=0.4$ [cm] and $x_e=2.2$ [cm], the velocity profiles match, as a fully developed flow state has been reached at both positions.}
\label{fig:cfd_flow_field_gnf}
\end{figure}

  \cref{fig:cfd_flow_field_gnf} shows the resulting flow field from a CFD simulation using $v^{\max}_{\text{inlet}}=0.03 \,[\frac{\text{m}}{\text{s}}]$ for a sample geometry. 
   The velocity $x$-component dominates the flow characteristic since the main flow direction follows the primary vessel axis.  
   The $y$-component becomes primarily non-zero within the stenotic region, where the flow adapts to the constricted artery. 
    Relative to the $x$-component, it is at least an order of magnitude smaller, depending on the severity of the stenoses.   
    In contrast to Newtonian fluids with constant viscosity, generalized Newtonian fluids might develop in straight channels a plug-like velocity profile rather than a parabolic one \cite{irgens_rheology}.
    The specific shape depends on the viscosity parameters and the local shear rate regime. 
    We examine arteries with a maximum diameter of $d_{\text{artery}}=1$ [mm] and maximum inflow velocities in the range $v^{\max}_{\text{inlet}} \in [0.02,0.6] \,[\frac{\text{m}}{\text{s}}]$,
    where the considered flows are laminar within a low Reynolds number regime. Importantly, we will restrict our study to geometries and flow conditions that admit a unique stationary solution.
    Our computational domains are designed with sufficient inlet length to reach fully developed flow profiles upstream of the stenoses at the cross-section $x=x_s$. Similarly, the extended outlet region allows the flow to reach a developed state at the downstream cross-section $x=x_e$. 
    Notably,  the stress-free boundary condition for generalized Newtonian fluids leads to a slight expansion of the flow at the outlet, however, the effects on the flow field further upstream from the outlet are negligible for the considered flows \cite{pachecognf}.

\paragraph{\textbf{ Generating Training Data: Numerical Methods and Software}}
 The training data for our surrogate model is generated using high-fidelity numerical simulation methods and software. Here, the coupled nonlinear PDE system \cref{eq:gnf_pde} is solved on a High-Performance Computing (HPC) cluster using a scalable domain decomposition approach within our in-house  finite element-based software library \texttt{FEDDLib} \cite{feddlib} in combination with the \texttt{FROSch} package \cite{Frosch} from the \texttt{Trilinos} library \cite{mayr2025trilinosenablingscientificcomputing}.
 In \cref{simulation}, tolerances of the linear and nonlinear solver and details of the HPC cluster are given.  For a more detailed description of the applied numerical methods and the software environment; see \cite{hochmuth, monolithic_preconditioner_lea}. 
  The continuous nonlinear problem \cref{eq:gnf_pde}  is first discretized using the finite element method (FEM) with P2-P1 elements, resulting in a discrete but still nonlinear system \cite{volker, article_he}.  The problem is then linearized using a combination of Picard's and Newton's method. Initially, the Picard method reduces the residual, followed by Newton's method to accelerate convergence. The resulting linear system is then solved efficiently with a domain decomposition preconditioning strategy, which exploits the parallel architecture of modern HPC clusters by decomposing the domain into overlapping subdomains. More precisely, the iterative solver GMRES (Generalized Minimal Residual) and an overlapping Schwarz preconditioner with a GDSW-type (Generalized Dryja-Smith-Widlund) coarse space implemented in the \texttt{Trilinos}-package \texttt{FROSch} is used \cite{monolithic_preconditioner_lea}. 
This  framework allows us to generate training data within a reasonable time frame.

\section{ The CNN-Schwarz-Flow Approach for Stationary Blood Flow Prediction}
\label{sec:new_algorithm_dd_cnn}

The CNN-Schwarz-Flow approach is based on an alternating Schwarz domain decomposition method to solve the  non-Newtonian flow problem on the original (global) domain $\Omega$. Here, $\Omega$ is decomposed into several overlapping subdomains and we need to solve smaller subproblems, related to the original non-Newtonian flow problem, on these subdomains. In order to define our global surrogate model, we will need to define such models for the subdomains. The strength of our approach is that we only need to train such a model on a single (reference) subdomain resulting in a universal subdomain solver (USDS). The inference of the surrogate model can then be used on all subdomains applying the USDS.
  Accordingly, we distinguish between two types of predictions: subdomain prediction and global prediction. For the subdomain prediction, the CNN architecture is trained and evaluated locally on pixel grids with fixed dimensions using data obtained from FEM simulations; see \cref{sec:single_image_prediction} for details. This results in the USDS. The global prediction based on the alternating Schwarz method extends to pixel grids of varying sizes by decomposing the grid into overlapping subdomains and combining multiple subdomain predictions by employing the locally trained 
  USDS on each subdomain; see \cref{sec:global_image_prediction} for details.

\subsection{Subdomain Prediction: CNN-based USDS}
\label{sec:single_image_prediction}

\begin{figure}[ht!]
\centering
    \resizebox{\textwidth}{!}{ 
\begin{tikzpicture}[scale=0.4]
    \node[anchor=south west] at (-0.95, -0.05) {\includegraphics[width=16cm]{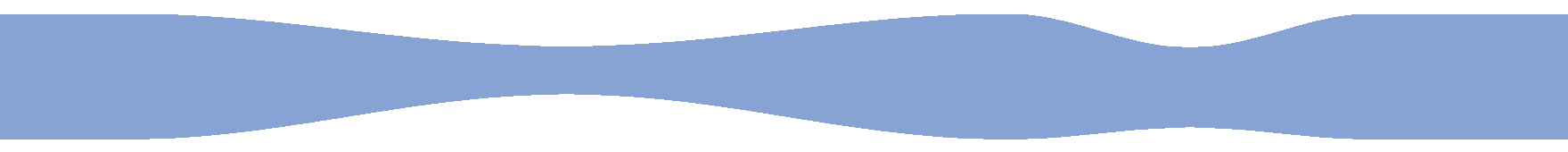}};
    
    \def\rectWidth{6.67}
    \def\rectHeight{3.3}
    \foreach \j  [count=\p from 0] in {-1, 0, 1, 2,3,4} {
            
                \begin{scope}
        \clip
            (\j * \rectWidth + 6, 0.5) rectangle 
            ({(\j + 1) * \rectWidth + 6}, \rectHeight + 0.5);
        \draw[step=0.2, thin, black!50] 
            (\j * \rectWidth + 6, 0.5) grid 
            ({(\j + 1) * \rectWidth + 6}, \rectHeight + 0.5);
    \end{scope}
    
     \draw[thick, black] 
            (\j * \rectWidth+6, 0.5) rectangle 
            ({(\j + 1) * \rectWidth + 6}, \rectHeight + 0.5);

        \node[above] at (\j * \rectWidth+6+3.5, 1.6) {\large $I_{\p}$ };
    }
\end{tikzpicture}
}
\vspace{-18pt}
\caption{Decomposition of an arterial geometry into six non-overlapping pixel grids \\$I_n, \,\, n=0,..., 5$. Each grid $I_n$ has a fixed dimension $ W\times H =256\times 128$  
pixels. }
\label{single_images}
\end{figure}

The input and output data for training the CNN are sampled from CFD simulations on rectangular pixel grids $I$ defined by \begin{equation}
I = \{(i,j) \in \mathbb{Z}^2 : 0 \leq i < W, 0 \leq j < H\},
\label{eq:image}
\end{equation}
with  a constant height of $H=128$ and width of $W=256$ pixels. Each element $(i,j) \in I$ represents a pixel location, with a total number of $|I| = W \times H$ pixels in the grid. The pixel grids are extracted from our two-dimensional computational domains.
These domains can be decomposed into pixel grids sequentially ordered along the length of the artery; see \cref{single_images}. 
This allows each pixel grid to cover the entire arterial diameter, with each geometry segment being completely embedded within the grid. 
Each pixel grid then serves as an independent training image for optimizing the weights of the CNN. In \cref{postprocessing}, the detailed procedure to extract the training data from the simulation results is described. 

\subsubsection{\textbf{Flow-Rate-Conserving CNN}}\label{sec_flow_rate_conserving_cnn}
\begin{figure}[ht!]
    \includegraphics[width=0.99\textwidth]{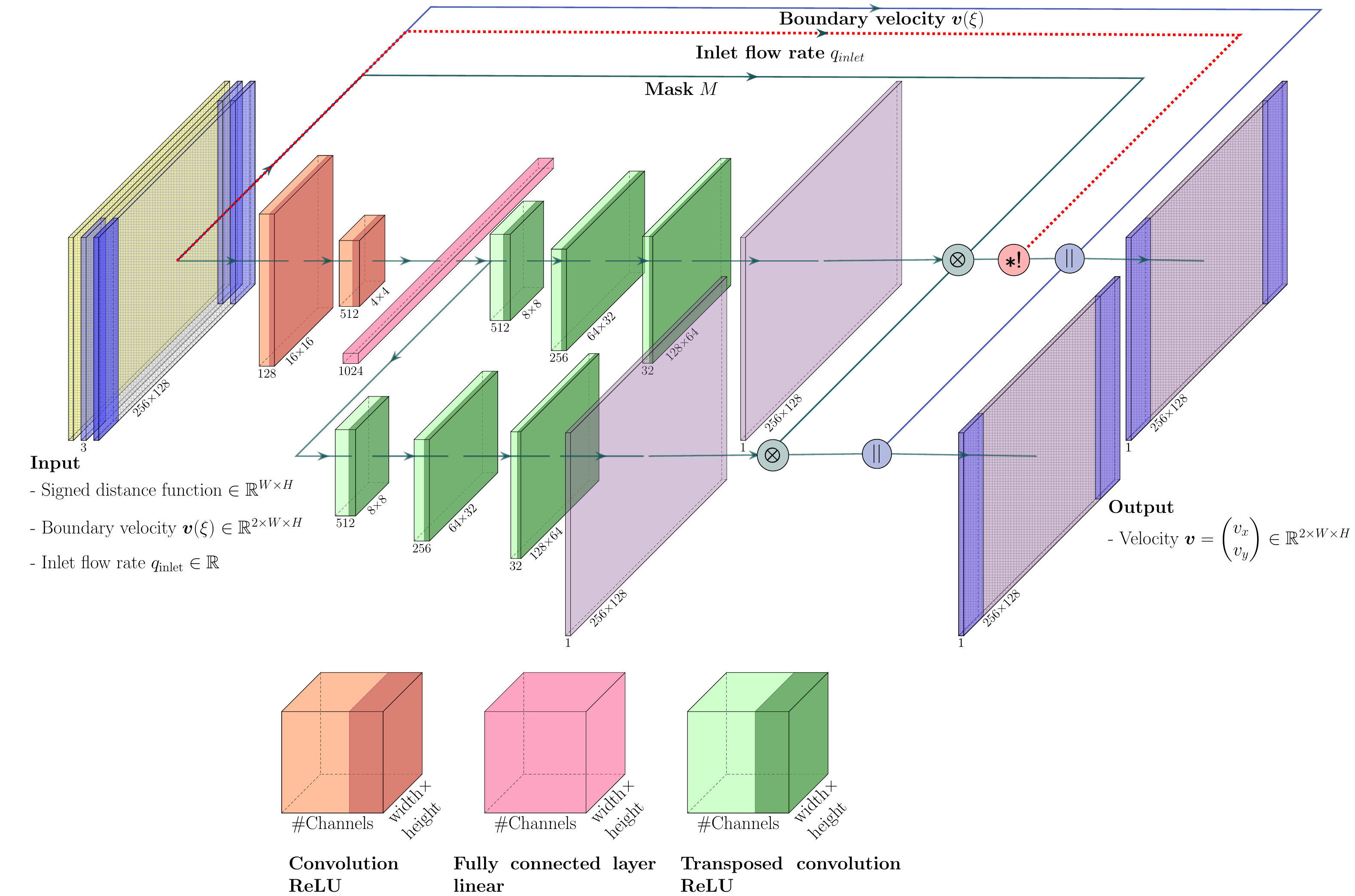}
        \vspace{-9pt}
\caption{Bottleneck CNN model for the prediction of the velocity components $v_x$ and $v_y$, adapted from \cite[fig. 9]{guo_cnn}.  
 The network employs a shared encoding structure followed by separated decoding branches. 
Color-coded operations include convolutions each followed by ReLU activation functions (darker colored), a fully connected layer for the latent space representation, and transposed convolutions followed by ReLU non-linearities in the decoders.
 The network takes three input channels of size $256 \times 128$; the SDF and the velocity components $\bm{v}(\xi)$ taken from a certain boundary region of the subdomain. Additionally, the inlet flow rate $q_{\text{inlet}}$ is required as scalar input, as it is used in the constraint layer (denoted by the operator $\textbf{*!}$) to enforce a constant flow rate in each cross-section using equation  \cref{constraint_equation}. 
 The network predictions undergo then in total three post-processing steps: First, the output is multiplied with the mask $M$ (denoted by the operator $\otimes$), then the constraint layer is applied, and finally, the input velocities  $\bm{v}(\xi)$ are explicitly set in the output (denoted by the operator $||$). The plot was generated using the tool \cite{plotnn}.}
\label{cnn_architecture}
\end{figure}

Our network architecture for the USDS is built upon the bottleneck CNN model introduced in \cite{guo_cnn};  see \cref{cnn_architecture}. 
Following the notation of \cite{eichinger_cnn}, the CNN is defined as a mapping from a three-channel input and an additional scalar input to a two-channel output.
Specifically, the network is trained to map image representations of an SDF and velocities from a boundary region of the pixel grid to an image representation of the velocity components, that is,
\begin{equation}
F_\text{CNN} : \mathbb{R}^{3 \times W \times H} \times \mathbb{R} \to \mathbb{R}^{2 \times W \times H},  \quad
\begin{bmatrix} 
SDF, 
v_x(\xi),
v_y(\xi) 
\end{bmatrix}^{\top} \times q_{\text{inlet}} \mapsto \bm{v} = \begin{bmatrix} v_x \\ v_y \end{bmatrix}.
\end{equation}
The inlet flow rate $q_{\text{inlet}}$ is provided as scalar input for the integration of a physics-aware constraint layer into the architecture; details are given in the following paragraphs.
The building blocks of the CNN architecture are multiple convolutional layers in a shared encoder, which extract spatial features from the image representation of the input data into a latent space. 
This information is then processed by two separate decoder paths with transposed convolutions to predict the velocity components independently. 
All layer characteristics, including kernel sizes and stride values, are maintained from the original work~\cite{guo_cnn}.

\paragraph{\textbf{Input: SDF}}
To predict the flow fields for different geometries, that is, varying stenosis strengths, an  SDF is provided as input. The SDF gives, in each pixel, information about its minimal distance to the boundary wall $\Gamma_w$. For two-dimensional grids, the boundary wall $\Gamma_w$ includes the upper and lower boundaries of the artery; see \cref{fig_combined} (left). More precisely, the SDF value in pixel $(i,j) \in I$ is defined as 
\begin{equation}
\label{eq:sdf}
\text{SDF}(i,j) = 
\begin{cases}
d((i,j), \Gamma_w) & \text{for } (i,j) \in \Omega \setminus \Gamma_w \\
0 & \text{for } (i,j) \in \Gamma_w \cup (I \setminus \Omega)
\end{cases},
\end{equation}
where $d( (i,j), \Gamma_w)$ is defined as shortest Euclidean distance from pixel $(i,j)$ to the boundary wall $\Gamma_w$. On $\Gamma_w$, as well as outside the geometry, the SDF is defined as zero. 
In \cref{sdf}, the 
method used for computing the SDF is mentioned.
\begin{figure}[h!]
\centering
\begin{minipage}{0.48\textwidth}
\centering
\begin{tikzpicture}[scale=0.5]
    \fill[black!20, opacity=0.5] (0, 0) rectangle (8, 4);
    
    \fill[blue!25, opacity=0.7] 
        plot[domain=0:8, samples=100] (\x, {3.2 + 0.2*(\x-4)^2/4}) -- 
        plot[domain=8:0, samples=100] (\x, {0.8 - 0.2*(\x-4)^2/4}) -- 
        cycle;

	\node[] at (2, 3.2) {...}; 
	\node[] at (2, 0.8) {...}; 

	\node[] at (6, 3.2) {...}; 
	\node[] at (6, 0.8) {...};

\node[rotate=90] at (0.5, 2) {...}; 
\node[rotate=90] at (7.5, 2) {...}; 

\foreach \x/\color in {0/red!30, 0.333/red!30, 0.666/gray!20, 1/gray!20} {
    \foreach \y in {3.75} {
        \draw[fill=\color] (\x, \y) rectangle (\x+0.333, \y+0.25);
    }
}

\foreach \x/\color in {0/blue!30, 0.333/blue!30, 0.666/red!30, 1/red!30} {
    \foreach \y in {3.5} {
        \draw[fill=\color] (\x, \y) rectangle (\x+0.333, \y+0.25);
    }
}

\foreach \x/\color in {0/blue!30, 0.333/blue!30, 0.666/blue!30, 1/blue!30} {
    \foreach \y in {3.25, 3.0, 2.75} {
        \draw[fill=\color] (\x, \y) rectangle (\x+0.333, \y+0.25);
    }
}

\foreach \x/\color in {6.666/blue!30, 7/blue!30, 7.333/blue!30, 7.666/blue!30} {
    \foreach \y in {1.25, 1.0, 0.75, 0.5} {
        \draw[fill=\color] (\x, \y) rectangle (\x+0.333, \y+0.25);
    }
}

\foreach \x/\color in {6.666/red!30, 7/red!30, 7.333/blue!30, 7.666/blue!30} {
    \foreach \y in {0.25} {
        \draw[fill=\color] (\x, \y) rectangle (\x+0.333, \y+0.25);
    }
}

\foreach \x/\color in {6.666/gray!20, 7/gray!20, 7.333/red!30, 7.666/red!30} {
    \foreach \y in {0} {
        \draw[fill=\color] (\x, \y) rectangle (\x+0.333, \y+0.25);
    }
}

\foreach \x/\color in {6.666/gray!20, 7/gray!20, 7.333/red!30, 7.666/red!30} {
    \foreach \y in {3.75} {
        \draw[fill=\color] (\x, \y) rectangle (\x+0.333, \y+0.25);
    }
}

\foreach \x/\color in {6.666/red!30, 7/red!30, 7.333/blue!30, 7.666/blue!30} {
    \foreach \y in {3.5} {
        \draw[fill=\color] (\x, \y) rectangle (\x+0.333, \y+0.25);
    }
}

\foreach \x/\color in {6.666/blue!30, 7/blue!30, 7.333/blue!30, 7.666/blue!30} {
    \foreach \y in {3.25, 3.0, 2.75} {
        \draw[fill=\color] (\x, \y) rectangle (\x+0.333, \y+0.25);
    }
}

\foreach \y in {0.5, 0.75, 1, 1.25} {
    \foreach \x in {0, 0.333, 0.666, 1} {
        \draw[fill=blue!30] (\x, \y) rectangle (\x+0.333, \y+0.25);
    }
}

\foreach \x/\color in {0/blue!30, 0.333/blue!30, 0.666/red!30, 1/red!30} {
    \draw[fill=\color] (\x, 0.25) rectangle (\x+0.333, 0.5);
}

\foreach \x/\color in {0/red!30, 0.333/red!30, 0.666/gray!20, 1/gray!20} {
    \draw[fill=\color] (\x, 0) rectangle (\x+0.333, 0.25);
}

    \draw[thick, red, dashed, opacity=1.0] plot[domain=0:8, samples=100] (\x, {3.2 + 0.2*(\x-4)^2/4});  
    \draw[thick, red, opacity=1.0, dashed] plot[domain=0:8, samples=100] (\x, {0.8 - 0.2*(\x-4)^2/4});

    \draw[decorate, decoration={brace, amplitude=10pt}, thick, black] (8, -0.05) -- (0, -0.05); 
    \node[below] at (4, -0.6) {\small $256$ Pixels}; 
  
    \draw[decorate, decoration={brace, amplitude=10pt}, thick, black] (-0.1, 0) -- (-0.1, 4); 
    \node[above, rotate=90] at (-0.9, 2) {\small $128$ Pixels}; 

    \draw[thick, black] (0, 0) rectangle (8, 4);
    \node[blue!80] at (4, 2) {$\Omega$}; 

    \draw[->,  red] (5.7, 2.4) -- (5, 3.1);
    \draw[->,  red] (5.7, 2.2) -- (5, 0.84);
    \node[red] at (6.3, 2.3) {$\Gamma_w$}; 
    \node[below] at (0.5, -0.5) {\vphantom{$\xi=2$ Pixels}};

\end{tikzpicture}
\end{minipage}
\hfill
\begin{minipage}{0.48\textwidth}
\centering
\begin{tikzpicture}[scale=0.5]

    \fill[black!20, opacity=0.5] (0, 0) rectangle (8, 4);

\foreach \x/\color in {0/blue!50, 0.333/blue!50, 0.666/gray!20, 1/gray!20} { 
    \foreach \y in {0, 0.25, ..., 1.3} { 
        \draw[fill=\color] (\x, \y) rectangle (\x+0.333, \y+0.25); 
    }
}
	\node[] at (2, 3.2) {...}; 
	\node[] at (2, 0.8) {...}; 

	\node[] at (6, 3.2) {...}; 
	\node[] at (6, 0.8) {...};

\node[rotate=90] at (0.5, 2) {...}; 

\foreach \x/\color in {0/blue!50, 0.333/blue!50, 0.666/gray!20, 1/gray!20} { 
    \foreach \y in {2.75, 3.0, ..., 3.9} { 
        \draw[fill=\color] (\x, \y) rectangle (\x+0.333, \y+0.25); 
    }
}

\foreach \x/\color in {6.666/gray!20, 7/gray!20, 7.333/blue!50, 7.666/blue!50} { 
    \foreach \y in {0, 0.25, ..., 1.3} { 
        \draw[fill=\color] (\x, \y) rectangle (\x+0.333, \y+0.25); 
    }
}

\node[rotate=90] at (7.5, 2) {...}; 

\foreach \x/\color in {6.666/gray!20, 7/gray!20, 7.333/blue!50, 7.666/blue!50} { 
    \foreach \y in {2.75, 3.0, ..., 3.9} { 
        \draw[fill=\color] (\x, \y) rectangle (\x+0.333, \y+0.25); 
    }
}

    \draw[thick, black] (0, 0) rectangle (8, 4); 
    \node at (4, 2) {$I$}; 

        \draw[<->, thick, blue] (0, -0.4) -- (0.6, -0.4);
            \node[below] at (0.5, -0.6) {\small $\xi=2$ Pixels};

		\draw[<->, thick, blue] (7.4, -0.4) -- (8, -0.4);
		\node[below] at (7.7, -0.6) {\small $\xi=2$ Pixels};
		
\end{tikzpicture}
\end{minipage}
\vspace{-7pt}
\caption{ Left: Visualization of an arterial segment within a pixel grid, showing the extracted fluid domain $\Omega$ bounded by the boundary wall $\Gamma_w$ (dashed lines). Right: Visualization of a pixel grid $I$ with colored pixel columns representing the set $B_{\xi=2}(I)$.}
\label{fig_combined}
\end{figure}

\paragraph{\textbf{Input: Velocity Components in Boundary Region}}
 While previous studies often only used an SDF or other geometry information as input, this work additionally incorporates velocity values from a boundary region of the pixel grid as input data. 
 This is a prerequisite for propagating and exchanging information when predicting stationary flows over multiple subdomains as described in \cref{sec:global_image_prediction}. Given our decomposition approach for arteries, subdomains can only overlap in the $x$-direction.
  Thus, the exchange of velocity information through overlaps of neighboring subdomains can occur exclusively through the left and right boundary regions of the subdomains. 
  Therefore, the input velocity components, that is, \(v_x(\xi)\) and \(v_y(\xi)\), are extracted here  only from the left and right boundary regions of the pixel grid. 
  These regions are parameterized by the input boundary width, denoted by $\xi$, which determines the number of pixel columns adjacent to the boundaries to be included.
 This is in contrast to classical boundary value problems, where the boundary condition is defined on a submanifold, providing the CNN with richer context about flow dynamics near the boundaries. 
  We will investigate in \cref{single_pred_performance}, if larger $\xi$ values improve the prediction accuracy. 
  To  formalize the regions contributing to the CNN input, 
  we define the set of pixels within distance $\xi$ from the left and right boundaries of $I$ as
\begin{equation}
B_{\xi}(I) = \{(i,j) \in I : 0 \leq i < \xi \text{ or } W-\xi \leq i < W, \, 0 \leq j < H\}.
\end{equation}
This set contains the first $\xi$ columns and the last $\xi$ columns of pixels in the grid.
For $\xi_{\min} = 1$, this set reduces to just the left and right pixel columns
$B_{1}(I)  = \{(i,j) \in I : i = 0 \text{ or } i = W-1, \, 0 \leq j < H\}$.
In \cref{fig_combined} (right), the set $B_{\xi=2}(I)$ is visualized for a grid $I$. The values for $(v_x(\xi), v_y(\xi))^{\top}  \in B_{\xi}(I) $ for the training are obtained by extracting these from the corresponding velocity results obtained from the FEM simulations. 
The remaining values are set to zero, that is, $(v_x(\xi), v_y(\xi))^{\top} = 0 \in I \setminus B_{\xi}(I)$.

\paragraph{\textbf{Input: Inlet Flow Rate for Physics-Aware Constraint Layer}}
\label{pyhsics}

Standard data-driven CNNs often generate predictions that violate physical principles such as mass conservation.
 To address this, two strategies are common: soft constraints,
  which add penalty terms based on physical laws to the loss function, and hard constraints, which directly incorporate physical laws into the model \cite{article_beucler_constraint}.
   While soft constraints usually fulfill the physical constraints in a least squares sense, hard constraints strictly enforce them.
   In this work, we adopt the latter approach and integrate a simple constraint layer into the CNN architecture to preserve the flow rate \cite{weather_cnn_downsampling}.
 
   For an incompressible, stationary flow in closed channels without sources or sinks, mass conservation implies that the flow rate $Q$ remains constant throughout the domain.
 For two-dimensional domains,
 $Q$ is defined as an integral of the velocity field $\bm{v}$ over a cross-sectional curve $C$, that is, 
\begin{equation}
Q = \int_{C} \bm{v} \cdot \bm{n} \, dC,
\end{equation}
where $\bm{n}$ is the unit outward normal vector to the curve \cite{kundu2024fluid}. For our two-dimensional geometries with horizontally aligned pixel grids, each cross-section is a vertical line perpendicular to the main flow direction (that is, $x$-direction) and the normal is directed in the $x$-direction.
Thus, only the $x$-component of the velocity $v_x(x,y)$ contributes to the flow rate, that is,
\begin{equation}
	Q(x) = \int v_x(x, y) \, dy,
\end{equation}
which is numerically approximated on the pixel grid $I$ by
\begin{equation}
    \label{eq:flowrate}
q(i) \approx \sum_{j=0}^{H-1} v_x(i, j) \cdot \Delta y, \quad 0 \leq i < W.
\end{equation}
Here, $i$ denotes the pixel column index  and $\Delta y = d_{\text{artery}} / H$ is the cell length corresponding to a pixel. 
Since, in our case, each arterial cross-section is completely embedded within a single subdomain, the predicted flow rates must be equal to the inlet flow rate $q_{\text{inlet}}$ for all $i$, where the value of $q_{\text{inlet}}$ can be determined from the parabolic inlet profile (equation \cref{parabolic}). 
However, without additional constraints, a CNN output $\tilde{\mathbf v}$ typically yields flow rates  $\tilde q(i)$ that do not strictly satisfy flow rate conservation, that is, $\tilde q(i) \neq q_{\text{inlet}}$.
Therefore, we define a scaling factor $C(i)=q_{\text{inlet}}/\tilde{q}(i)$ for each pixel column $i$ to capture these deviations.
We then multiply the predicted velocities $\tilde{v}_x(i,j)$ by these scaling factors $C(i)$  to enforce flow rate conservation at each cross-section, that is,
\begin{equation}
\tilde{v}^{\text{*}}_x(i,j) = C(i) \cdot \tilde{v}_x(i,j), \quad 0 \leq i < W, \, 0 \leq j < H.
\label{constraint_equation}
\end{equation}
 This ensures that the corrected flow rates strictly fulfill $\tilde{q}^{\text{*}}(i)=q_{\text{inlet}}$ for all $i$.
 Importantly, the loss function for training the CNN is evaluated based on the final prediction $\bm{v}$, including the scaling operation.
 Since the constraint layer depends on the network parameters through $\tilde{\bm{v}}$ and is part of the computational graph, 
 the gradient $\partial \tilde{\bm{v}}^{\text{*}} / \partial \tilde{\bm{v}}$ is required during backpropagation such that the constraint layer influences the parameter updates.
As the operation \cref{constraint_equation} is differentiable, the gradients can be automatically computed by \texttt{PyTorch}'s autograd framework.
We want the constraint layer to make only minor quantitative corrections to $\tilde{\bm{v}}$, which means that the CNN should already generate accurate predictions before scaling.
Thus, to prevent the constraint layer from potentially masking poor predictions $\tilde{\bm{v}}$ during training by computing, for example, extremely large or small scaling factors $C(i)$, 
 a penalty term is added to the loss function that penalizes large deviations of $C(i)$ from $1$, essentially imposing $\tilde{v}_x\approx\tilde{v}^{\text{*}}_x $; see \cref{training}.
 During inference, the values of $C(i)$ are monitored, and if any $C(i)<0$, a warning is triggered and the prediction is marked as potentially failed, 
 since negative scaling factor indicate an incorrectly predicted net flow direction and would reverse the flow direction.

  Notably, our approach only constrains the $x$-component of the velocity.
  A more general approach would be to incorporate mass conservation in terms of the divergence\allowbreak-free condition \cref{eq:divergendefree} into the model.
   In previous works, this has been achieved by predicting a vector potential and computing its curl to obtain a divergence\allowbreak-free velocity field \cite{Kim2019DeepFluids, PhysRevFluids.8.014604}. 
   For approximating the spatial derivatives in the curl operator the authors use non-trainable convolutional layers whose kernels correspond to finite-difference stencils with fixed weights. 
   Given our problem assumptions, we have decided to first examine a simpler approach that, as described, uses a linear scaling operation to ensure flow rate conservation without approximating spatial derivatives.
   Investigating a divergence-free CNN architecture as USDS within the proposed Schwarz framework is left for future work.

\paragraph{\textbf{Output: Velocity components}}

 The outputs of the CNN are the velocity field components $\bm{v} =  (v_x,  v_y )^{\top}
$. As in \cite{guo_cnn}, the binary mask $M:  I \rightarrow \{0,1\}$ is  defined using the SDF values
 to guarantee that the predicted velocities fulfill the no-slip boundary condition at the boundary of the arterial wall and are zero outside of it, 
 that is, 
\begin{equation}
\label{eq:mask}
M(i,j) =
\begin{cases}
1 & \text{if } \text{SDF}(i,j) > 0 \\
0 & \text{if } \text{SDF}(i,j) = 0.
\end{cases}
\end{equation} 
The predicted velocity components are multiplied pixel-by-pixel with this mask as a postprocessing step to the network output \cite{guo_cnn}. 
This is followed by applying the constraint layer \cref{constraint_equation} to the $x$-component of the predicted velocity field to enforce flow rate conservation.
Finally, we explicitly prescribe the provided input velocities $\bm{v}(\xi)$  in the network output,
extending the concept of Dirichlet boundary conditions to the entire boundary region $B_{\xi}(I)$, that is,
\begin{align}
\bm{v} (i,j) &= \bm{v}(\xi) (i,j) \quad \text{for } (i,j) \in B_{\xi}(I).
\end{align}
This is motivated by the original alternating Schwarz method, described in the following  \cref{sec:global_image_prediction},
in which Dirichlet boundary conditions are imposed in each subdomain.

\paragraph{\textbf{Remark: Data-driven CNN}}
To evaluate the impact of a physics-constrained design on prediction accuracy and convergence, 
we will compare the flow\allowbreak-rate\allowbreak-conser\-ving USDS against a purely data-driven variant. 
The data-driven CNN has the same encoder-decoder architecture (\cref{cnn_architecture}), but does not include the constraint layer and thus does not require the inlet flow rate $q_{\text{inlet}}$ as input.

\subsection{Global Prediction: Alternating Schwarz Method}
\label{sec:global_image_prediction}

In order to make a prediction for longer and more complicated arterial geometries without the necessity to train  
a new CNN, a Schwarz 
domain decomposition strategy is applied.
The alternating Schwarz method was originally introduced in 1870 by Hermann Amandus Schwarz for the solution of boundary value problems based on PDEs; see the monography \cite{book_domaindecomposition} for further references. 
Schwarz's key idea was to decompose a global domain into overlapping subdomains, define local problems in each subdomain, solve them iteratively, and exchange boundary information across the overlap within each iteration. 
For elliptic Dirichlet PDE problems, it can be proven
that the iterative procedure converges to the global solution, if the subdomain problems are solved exactly; see, for example, \cite{book_domaindecomposition} for further details and references. 
In this work, we replace the local subdomain solvers with the locally trained CNN surrogate models - or universal subdomain solver (USDS) - introduced 
in \cref{sec:single_image_prediction},  which can be considered as inexact solvers. Let us note again that the USDS is only trained on a single reference subdomain and then used for all subdomains in the alternating Schwarz method; no additional training for the subdomains used in the Schwarz method has to be carried out. Since there is no proof of convergence for the alternating Schwarz algorithm using inexact subdomain solvers, empirical tests must be performed to investigate the convergence behavior of the method.

\begin{figure}[h]
\centering
    \resizebox{\textwidth}{!}{ 
\begin{tikzpicture}[scale=0.32]

    \node[anchor=south west, opacity=0.5] at (-0.42, -0.7) {\includegraphics[width=11.1cm,height=1.5cm]{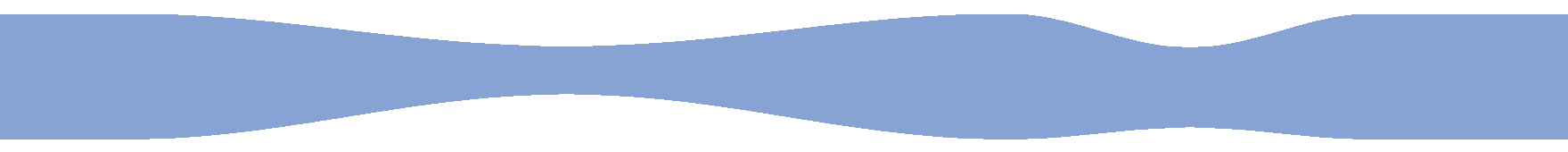}};

  \node[below] at (7.5, -0.5) {\small $\delta=4$ Pixels};
    \draw[<->, thick, black] (6.666, -0.4) -- (8, -0.4);
    
    \node[below] at (14, -0.5) {\small $\delta=4$ Pixels};
    \draw[<->, thick, black] (13.3333, -0.4) -- (14.6666, -0.4);
 
    \node[below] at (20.666, -0.5) {\small $\delta=4$ Pixels};
    \draw[<->, thick, black] (20, -0.4) -- (21.333, -0.4);

    \node[below] at (27.333, -0.5) {\small $\delta=4$ Pixels};
    \draw[<->, thick, black] (26.666, -0.4) -- (28, -0.4);
 
    \fill[red!20, opacity=0.5] (0, 0) rectangle (8, 4);
    \fill[black!20, opacity=0.5] (6.666, 0) rectangle (14.666, 4);
    \fill[red!20, opacity=0.5] (13.333, 0) rectangle (21.333, 4);
    \fill[black!20, opacity=0.5] (20, 0) rectangle (28, 4);
    \fill[red!20, opacity=0.5] (26.666, 0) rectangle (34.666, 4);
 
    \foreach \x in {6.666, 7, 7.333, 7.666} { 
        \foreach \y in {0, 0.25, ..., 1.3} { 
            \draw[] (\x, \y) rectangle (\x+0.333, \y+0.25); 
        }
    }
    \node[rotate=90] at (7.4, 2) {...}; 
    \foreach \x in {6.666, 7, 7.333, 7.666} { 
        \foreach \y in {2.75, 3.0, ..., 3.9} { 
            \draw[] (\x, \y) rectangle (\x+0.333, \y+0.25); 
        }
    }
 
    \foreach \x in {13.333, 13.666, 14, 14.333} {
        \foreach \y in {0, 0.25, ..., 1.3} { 
            \draw[] (\x, \y) rectangle (\x+0.333, \y+0.25); 
        }
    }
    \node[rotate=90] at (14.1, 2) {...}; 
    \foreach \x in {13.333, 13.666, 14, 14.333} { 
        \foreach \y in {2.75, 3.0, ..., 3.9} { 
            \draw[] (\x, \y) rectangle (\x+0.333, \y+0.25); 
        }
    }
 
    \foreach \x in {20, 20.333, 20.666, 21} {
        \foreach \y in {0, 0.25, ..., 1.3} { 
            \draw[] (\x, \y) rectangle (\x+0.333, \y+0.25); 
        }
    }
    \node[rotate=90] at (20.7, 2) {...}; 
    \foreach \x in {20, 20.333, 20.666, 21} { 
        \foreach \y in {2.75, 3.0, ..., 3.9} { 
            \draw[] (\x, \y) rectangle (\x+0.333, \y+0.25); 
        }
    }

    \foreach \x in {26.666, 27, 27.333, 27.666} {
        \foreach \y in {0, 0.25, ..., 1.3} { 
            \draw[] (\x, \y) rectangle (\x+0.333, \y+0.25); 
        }
    }
    \node[rotate=90] at (27.4, 2) {...}; 
    \foreach \x in {26.666, 27, 27.333, 27.666} { 
        \foreach \y in {2.75, 3.0, ..., 3.9} { 
            \draw[] (\x, \y) rectangle (\x+0.333, \y+0.25); 
        }
    }
 
    \draw[thick, red] (0, 0) rectangle (8, 4); 
    \node at (4, 2) {$I_0'$};
 
    \draw[thick, black] (6.666, 0) rectangle (14.6666, 4); 
    \node at (11, 2) {$I_1'$};
 
    \draw[thick, red] (13.333, 0) rectangle (21.333, 4); 
    \node at (18.0, 2) {$I_2'$};
 
    \draw[thick, black] (20, 0) rectangle (28, 4); 
    \node at (24.5, 2) {$I_3'$};

    \draw[thick, red] (26.666, 0) rectangle (34.666, 4); 
    \node at (31.2, 2) {$I_4'$};
 
    \draw[decorate, decoration={brace, amplitude=10pt}, thick, red] (0, 4.2) -- (8, 4.2); 
    \node[above] at (4, 5.2) {\small $256$ Pixels};
    
    \draw[decorate, decoration={brace, amplitude=10pt}, thick, black] (6.666, 4.2) -- (14.666, 4.2); 
    \node[above] at (11, 5.2) {\small $256$ Pixels}; 
    
    \draw[decorate, decoration={brace, amplitude=10pt}, thick, red] (13.333, 4.2) -- (21.333, 4.2); 
    \node[above] at (18, 5.2) {\small $256$ Pixels}; 
 
    \draw[decorate, decoration={brace, amplitude=10pt}, thick, black] (20, 4.2) -- (28, 4.2); 
    \node[above] at (24.5, 5.2) {\small $256$ Pixels}; 

    \draw[decorate, decoration={brace, amplitude=10pt}, thick, red] (26.666, 4.2) -- (34.666, 4.2); 
    \node[above] at (31.2, 5.2) {\small $256$ Pixels}; 
    
    \draw[decorate, decoration={brace, amplitude=10pt}, thick, red] (-0.1, 0) -- (-0.1, 4); 
    \node[above, rotate=90] at (-1.2, 2) {\small $128$ Pixels};
\end{tikzpicture}
    }
\vspace{-10pt}
\caption{Visualization of five overlapping subdomains $I_i', \, i=0,\ldots,4$ with each $256 \times 128$ pixels and an overlap width $\delta=4$, covering a global pixel grid with  $1264 \times 128$ pixels.}
\label{overlap}
\end{figure}

\paragraph{\textbf{Domain Decomposition}} 
We constrain ourselves to arterial geometries that can be embedded in a global pixel grid with $W_{\text{global}} \times  H $ pixels. 
For our application, the number of pixels in height corresponds to that for a local subdomain $H=128$, 
while the number of pixels in width $W_{\text{global}}$ and thus the length of the embedded arterial geometry can vary.  
 The global grid is then decomposed into $N$ overlapping subdomains $\{I_n'\}_{n=0}^N$, which overlap in the $x$-direction. 
 Here, the overlap width $\delta$ is defined as the number of pixel columns  shared between two adjacent subdomains. 
  \cref{overlap} shows the decomposition of a sample global pixel grid using five subdomains and an overlap width $\delta = 4$.
The feasible values of the overlap width depend on $\xi$, with $\delta_{\min} = 2\xi$ and in our case $\delta_{\max} = W - \xi$, as will become clear in the next section. 
Since $\delta_{\min} < \delta_{\max}$, the theoretical limit for $\xi$ is in our case $\xi < W/3 < 86$.

\paragraph{\textbf{CNN-Schwarz-Flow Prediction Approach}}

For  a new global prediction, an alternating Schwarz approach is employed, where the local subdomain solver is replaced by the evaluation of the USDS, our 
locally trained CNN surrogate model.  
Specifically, a red-black decomposition scheme 
is used in which the subdomains are partitioned into even-numbered (red) and odd-numbered (black) subdomains \cite{book_domaindecomposition}; see \cref{redblack_algo} for details. 
 In the red phase, the initially available information is used in the even-numbered subdomains as input to predict the interior flow fields. Subsequently, in the black phase, updated velocity information from neighboring subdomains is used to predict the flow fields within the black subdomains. Here, the overlap width constraints $\delta_{\min}$ and $\delta_{\max} $ become important since the overlap has to be chosen in relation to the input width $\xi$ such that it is guaranteed that all pixels are updated in each Schwarz iteration.
An iteration of the Schwarz method is completed with each red-black cycle. Through successive iterations, information about boundary conditions is progressively propagated from the global boundaries to the interior, and interior subdomain solutions are updated accordingly. Based on a specified stopping criterion and threshold $\epsilon$,
  this iterative loop should proceed until the predicted solution converges to an approximation of the global solution that satisfies the specified global boundary conditions.
\begin{algorithm}[]
 \caption{\textbf{CNN-Schwarz-Flow Prediction Approach with Red-Black Partitioning Using Locally Trained CNNs as Inexact Subdomain Solvers}}
\begin{algorithmic}[1]
\STATE{Compute signed distance function (SDF)}
\STATE{Decompose global domain into overlapping subdomains $I_n'$ with $\delta_{\min} \geq 2\xi  $}
\STATE{Boundary conditions: Prescribe velocities in global boundary regions }
\STATE{Initialization: Set for interior subdomains initial input velocities} in $B_{\xi}(I_n')$
\WHILE{Stopping criterion  $> \epsilon$}
  \IF{even $n$ subdomains {\color{red}(\textbf{Red phase})}}
       \STATE{ USDS: $F_{CNN} \left( \begin{bmatrix} SDF, v^k_x(\xi), v^k_y(\xi) \end{bmatrix}^{\top}_n \times q_{\text{inlet}} \right)  =  \begin{bmatrix} v^{k+\frac{1}{2}}_x \\ v^{k+\frac{1}{2}}_y \end{bmatrix}_n $} 
    \ELSIF{odd $n$ subdomains (\textbf{Black phase})}
        \STATE{ USDS: $F_{CNN}$ $\left(\begin{bmatrix} SDF, \,\, v^{k+\frac{1}{2}}_x(\xi), \,\,  v^{k+\frac{1}{2}}_y(\xi),\end{bmatrix}^{\top}_n \times q_{\text{inlet}}  \right)  =  \begin{bmatrix} v^{k+1}_x \\ v^{k+1}_y \end{bmatrix}_n $} 
    \ENDIF
    \STATE{k = k+1}
\ENDWHILE
\end{algorithmic}
\label{redblack_algo}
\end{algorithm}

\paragraph{\textbf{ Remark: Stopping Criterion}}
We define a stopping criterion for \cref{redblack_algo}  based on the update of the global prediction $\delta \bm{v} = \mathbf{v}^{k+1} - \mathbf{v}^k$. Several approaches can be employed, such as a criterion based on the $\ell_2$-norm $\|\delta \bm{v}\|_2$. However, the $\ell_2$-norm scales with the number of pixels, making it challenging to establish a consistent threshold $\epsilon$ across varying domain sizes.  Therefore, a stopping criterion based on the maximum absolute difference between consecutive Schwarz iterations is used instead, that is, 
\begin{equation}
    \label{eq:max:abs_error}
\text{AbsErr}_{\max} = \max_{i,j} |\bm{v}^{k+1}_{i,j} - \bm{v}^k_{i,j}|, \quad 0 \leq i < W_{\text{global}}, \, 0 \leq j < H,
\end{equation}
which is independent of the total pixel count.  In particular, by monitoring the maximum point-wise change, the iteration process stops if the largest update in any pixel falls below the specified threshold value $\epsilon$. Notably, using the USDS, that is, a CNN, as inexact subdomain solver offers no mathematical guarantee of convergence to the global solution, meaning that even if the stopping criterion reached $\epsilon$, the solution could have converged to an incorrect solution.
Therefore, in practice, it is essential to define an admissible parameter space regarding the geometric and flow configurations used during training of the USDS. Empirical tests must be performed to validate the convergence behavior of the CNN-Schwarz-Flow approach across these configurations.
In particular, it must be verified that the local prediction accuracy of the USDS is sufficient for the Schwarz iteration to converge, and that 
the training data is balanced so that all relevant flow features are adequately represented in the training data.
\paragraph{\textbf{Remark: Initialization of Input Velocities for Interior Subdomains}}
In practice, information on the velocity field is only known at cross-sections in the inlet and outlet region, where the flow in the straight channel segments has reached a developed state. 
 Therefore,  in our case, velocities are fixed in the boundary region of the first and last subdomain, and the information about the boundary conditions must then be propagated iteratively to all interior subdomains.
  To provide the USDS with a better initial guess than zeros, we set initial velocity values in the boundary regions $B_{\xi}(I_{n}')$ of the interior subdomains; see \cref{fig:initialization}.
   In this work, we initialize the velocity $x$-component with a parabolic profile, as it provides for the considered channel-like flows a reasonable starting value.
    The profile is computed based on a general parabolic equation, which takes into account the height of the arterial geometry in the considered cross-sections and which maintains the inlet flow rate $q_{\text{inlet}}$. 
    The $y$-component is initialized with zeros.
Alternative approaches could be applied, such as using constant values. However, poor initialization may lead to divergent predictions.  

We emphasize that the initialization is particularly important when using the flow-rate-conserving CNN with the global inlet flow rate $q_{\text{inlet}}$ as input.
As we do not train the CNN with zero velocity boundary conditions, initializing $v_x(\xi)=v_y(\xi)=0$ leads to intermediate predictions $\tilde{\bm{v}}$ with flow rates near zero, which causes the constraint layer to compute large scaling factors $C(i)$. 
These factors are then applied regardless of the prediction quality, either scaling a qualitatively correct prediction to the right magnitude or amplifying noisy predictions.
This could potentially provide neighboring subdomains with highly perturbed input velocity fields, which could lead to divergent predictions.
The proposed parabolic initialization prevents this issue.

 \begin{figure}[ht!]
\centering
    \resizebox{\textwidth}{!}{ 
\begin{tikzpicture}
    \node[anchor=south west,inner sep=0] (image) at (0,0) 
    {\includegraphics[scale=1.0]{figures/InitialBC/parabolic_init.png}};

    \def\xStart{17}   
    \def\yMin{53}    
    \def\yMax{79}    
    \def\width{52.5}   
    \def\overlap{12.3} 
    \foreach \i in {0,...,10} {
        \pgfmathsetmacro\xLeft{\xStart + \i * (\width - \overlap)}
        \pgfmathsetmacro\xRight{\xLeft + \width}
        
        \ifnum \i=0
            \def\rectColor{red}
             \draw[\rectColor, opacity=1.0, line width=20pt, dash pattern=on 100pt off 130pt] 
                    (\xLeft, \yMin) rectangle (\xRight, \yMax);

        \else
            \ifodd\i
                \def\rectColor{black!90!white}
                \draw[gray!30!white, opacity=1.0, line width=20pt, dash pattern=on 140pt off 160pt] 
                    (\xLeft, \yMin) rectangle (\xRight, \yMax);
            \else
                \def\rectColor{red}
                \draw[\rectColor, opacity=1.0, line width=20pt, dash pattern=on 100pt off 120pt] 
                    (\xLeft, \yMin) rectangle (\xRight, \yMax);
            \fi
        \fi

            \draw[decorate, decoration={brace, amplitude=140pt}, \rectColor, line width=20pt] 
           (\xLeft, \yMax+5) -- (\xRight, \yMax+5);
        
        \pgfmathsetmacro\bracketMid{(\xLeft + \xRight) / 2}
        \node[\rectColor] at (\bracketMid, \yMax+18) {\scalebox{32}{{ $I'_{\i}$}}};
    }

            \draw[decorate, decoration={brace, amplitude=100pt}, line width=20pt] 
           (149.5, 49) -- (137, 49);
   \node[below] at (142, 45) {\scalebox{22}{\large Overlap $\delta=60$}}; 

            \draw[decorate, decoration={brace, amplitude=20pt}, line width=20pt] 
           (20, 49) -- (16.5, 49);
             \node[below] at (20, 45) {\scalebox{22}{\large $\xi=20$}};

           \node[above] at (50, 110) {\scalebox{22}{\large Developed flow}}; 
           \node[below] at (440, 35) {\scalebox{22}{\large Developed flow}}; 

            \draw[->, line width=20pt, >=stealth] (24, 107) -- (20,81) ;

            \draw[->, line width=20pt, >=stealth] (465, 35) -- (469,48) ;

   \node[above] at (250, 40) {\scalebox{22}{\large \textbf{Parabolic initialization}}};

   \node[above] at (250, 22) {\scalebox{22}{\large Velocity $x$-component $[\frac{\text{m}}{\text{s}}]$}}; 

\end{tikzpicture}
}
\vspace{-15pt}
\caption{Illustration of the velocity field initialization for the first iteration of \cref{redblack_algo}.
The arterial geometry is embedded in a global pixel grid, which is decomposed into $N=11$ overlapping subdomains with $\xi=20$ and $\delta=60$.
In the first and last subdomain, the known developed flow profiles in the straight inlet and outlet segments are prescribed and remain fixed throughout the iterations.
For the interior subdomains, $v_x$ is initialized with a parabolic profile computed based on the local cross-section height and the known inlet flow rate, while $v_y$ is set to zero.}

\label{fig:initialization}
\end{figure}

\paragraph{\textbf{Remark: Scalability of CNN-Schwarz-Flow Approach}}

If a purely data-driven CNN is used in \cref{redblack_algo} as the subdomain solver, the approach corresponds to a one-level Schwarz method, since interior subdomains lack global information, 
and information on the global boundary conditions must be iteratively propagated inwards.
This makes the convergence of the method dependent on the number of subdomains. If $N$ is 
the number of subdomains and the first subdomain is labeled red, the approach requires $\lceil \frac{N}{2} \rceil$ red-black iterations
 to propagate the inflow information once through the entire domain.
 However, if, in \cref{redblack_algo}, the flow-rate-conserving CNN is used as the subdomain solver instead,
 each  interior subdomain already receives global information through the inlet flow rate $q_{\text{inlet}}$ in the constraint layer. 
  In \cref{scalability}, we investigate in initial tests how this affects the convergence with an increasing number of subdomains.

\section{Experimental Setup}
\label{sec:setup}

Following, the entire implemented pipeline for generating the training data and training the CNN-based USDS model is presented. The pipeline is divided into three distinct steps, each of which is explained in detail. 
\vspace{5pt}
   \begin{enumerate}
            \item \textbf{ Generation of training data} (\cref{generation_train})
            \begin{itemize}
            \item \textbf{Mesh generation}: \texttt{Pygmsh} package (\cref{geometry_fluidmodel})
            
            \item  \textbf{FEM simulation}: \texttt{FEDDLib}  (\cref{simulation})
            
            \item  \textbf{Grid sampling}: \texttt{PyVTK} package (\cref{postprocessing}) \begin{itemize}
            \item Computation of SDF
            \item Sampling of simulation data on structured pixel grids
            \item Extraction of subdomains for training
            \end{itemize}
            \end{itemize}
            \item \textbf{Training details} (\cref{training})
            \item \textbf{Application and analysis of results} (\cref{sec:results})
         \end{enumerate}

\subsection{Generation of Training Data}
\label{generation_train}

Generating a sufficient amount of high-quality training data, which covers a wide range of possible geometry and inflow configurations, is the most critical  and time-consuming part. Therefore, automation strategies were implemented to facilitate the computation of more training data.

\subsubsection{Mesh Generation}
\label{geometry_fluidmodel}

The pipeline starts with an automated mesh generation strategy using a Python script based on the \texttt{Pygmsh} package \cite{schlomer_pygmsh_2022}, followed by conversion of the meshes from the mesh format of \texttt{Gmsh} (\texttt{.msh}) to the \textit{MEDIT Inria} mesh format (\texttt{.mesh})  \cite{frey2001medit} using \texttt{MATLAB} \cite{MATLAB}. Each geometry consists of three main parts:  an inlet region, a stenotic region, and an outlet region.
 All geometries have a constant maximum diameter of $d_{\text{artery}}=1$ [mm] and a total length of $2.4$ [cm]. The stenotic region has a constant length of $1$ [cm] and can consist of up to three symmetric or asymmetric stenotic segments, each with a randomly assigned length and strength. In particular, the following parameters can vary between individual geometries
 \begin{itemize}
\item number of stenotic segments $n_{\text{segment}} \in \{1, 2, 3\}$,
\item relative length of $i$-th stenotic segment $l_{\text{segment, i}} \in [0.4, 1.5]$,
\item stenosis strengths for upper and lower walls $f_{\text{lower,i}}, f_{\text{upper,i}} \in [0.05, 0.75]$
\end{itemize}
with $i=1,\dots, n_{\text{segment}}$. These values are set randomly using the \texttt{random} module from Python. The relative lengths of the segments are normalized so that their sum matches the fixed length of the stenotic region. For each stenotic segment $i = 1, \dots, n_{\text{segment}}$ the lower and upper walls are defined by
\begin{subequations}
\label{stenosis}
\begin{align}
y_{\text{lower},i}(x) &= \frac{f_{\text{lower},i}}{2} \cdot \frac{d_{\text{artery}}}{2} \left(1 + \cos\!\left(\frac{2\pi(x - x_{0,i})}{l_{\text{segment},i}}\right)\right), \\
y_{\text{upper},i}(x) &= d_{\text{artery}} - \frac{f_{\text{upper},i}}{2} \cdot \frac{d_{\text{artery}}}{2} \left(1 + \cos\!\left(\frac{2\pi(x - x_{0,i})}{l_{\text{segment},i}}\right)\right),
\end{align}
\end{subequations}
with $x_{0,i}$ being the reference $x$-coordinate where the segment $i$ starts.  In total,  $770$ different meshes were constructed with each mesh consisting of approximately $1.5 \cdot 10^5$ triangular elements.
To characterize the geometric configurations in the dataset, we define two dimensionless parameters for each geometry; the normalized maximal severity of stenosis $S_{\text{norm}}$ and the normalized maximal stenosis gradient $G_{\text{norm}}$.
The parameter $S_{\text{norm}}$ captures the maximum severity of stenosis across all segments of a geometry, that is,
\begin{equation}
    S_{\text{norm}} = \max_{i \in \{1,\dots,n_{\text{segment}}\}} \frac{f_{\text{lower},i} + f_{\text{upper},i}}{f_{\text{lower}}^{\text{max}} + f_{\text{upper}}^{\text{max}}}  \in (0, 1],
\end{equation}
with $f_{\text{lower}}^{\text{max}} = f_{\text{upper}}^{\text{max}} = 0.75$.
The parameter $G_{\text{norm}}$ additionally accounts for the segment length, measuring how steeply the vessel narrows, that is,
\begin{equation}
    G_{\text{norm}} = \max_{i \in \{1,\dots,n_{\text{segment}}\}} S_{\text{norm}}\cdot\frac{l_{\text{segment},i}}{l^{\min}_{\text{segment}}}  \in (0, 1],
\end{equation}
where $l^{\min}_{\text{segment}}=0.4$ is the minimum possible relative segment length.
Higher values of $G_{\text{norm}}$ indicate steeper local constrictions, which are associated with more complex flow features such as stronger recirculation zones for higher velocities. 

\subsubsection{FEM Simulation}
\label{simulation}

Next, the meshes are used to conduct $770$ flow simulations, each with a randomly prescribed inlet velocity $v^{\max}_{\text{inlet}} \in [0.02, 0.6] \, [\frac{\text{m}}{\text{s}}]$ to capture different inflow regimes. 
Analogously to the normalized geometry parameters, we define a normalized flow rate parameter $Q_{\text{norm}}$ to characterize the inflow conditions, that is, 
\begin{equation}
    Q_{\text{norm}} = v^{\text{max}}_{\text{inlet}}/v^{\text{max}}_{\text{ref}} \, \in (0, 1], \quad \text{with } v^{\text{max}}_{\text{ref}} = 0.6 \, [\text{m/s}].
\end{equation}

The simulations are performed within our in-house software library \texttt{FEDDLib} and
 solve the system \cref{eq:gnf_pde} on a P$2$-P$1$ finite element grid for velocity and pressure degrees of freedom. The tolerance for the nonlinear and linear solver is set to $10^{-6}$, using a relative residual norm in both cases.
Compared to our initial work \cite{preprint}, the nonlinear residual tolerance was reduced from $10^{-4}$ to $10^{-6}$ to improve the solution accuracy, particularly for $v_y$.
 We observed that the higher tolerance led to numerical noise affecting the prediction accuracy of the CNN, especially for the purely data-driven CNN, which must learn the physics exclusively from the data.
The new simulations were performed on the Fritz supercomputer operated by the Friedrich-Alexander-Universit\"at Erlangen-N\"urnberg. The simulation ran on one node with two
Intel Xeon Platinum $8360$Y (Ice Lake) processors using $36$ MPI ranks, with each simulation requiring approximately $1$ minute.
 The resulting velocity fields are stored in XDMF/HDF5 files.

\subsubsection{Data Heatmaps: Geometry and Inflow Parameter Spaces}
\label{sec:dataheatmaps}

\begin{figure}[htbp]
    \begin{minipage}[t]{0.48\textwidth}
        \centering
        \includegraphics[width=\textwidth]{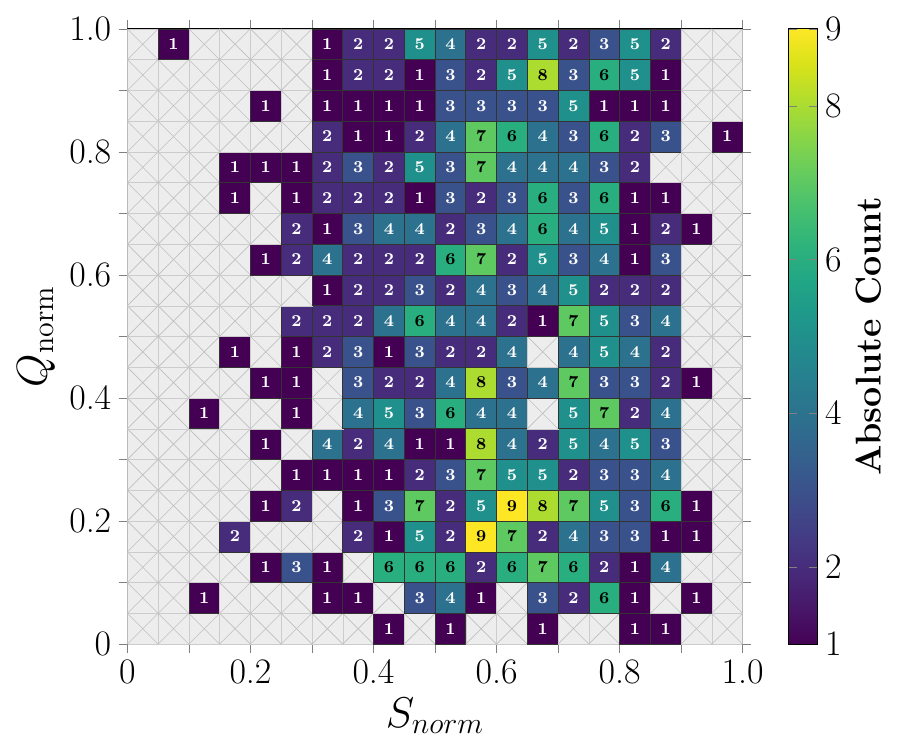}
        \par\smallskip
        \footnotesize\textit{(a) $S_{\text{norm}}$-$Q_{\text{norm}}$ parameter space.}
    \end{minipage}
    \hfill
    \begin{minipage}[t]{0.49\textwidth}
        \centering
        \includegraphics[width=\textwidth]{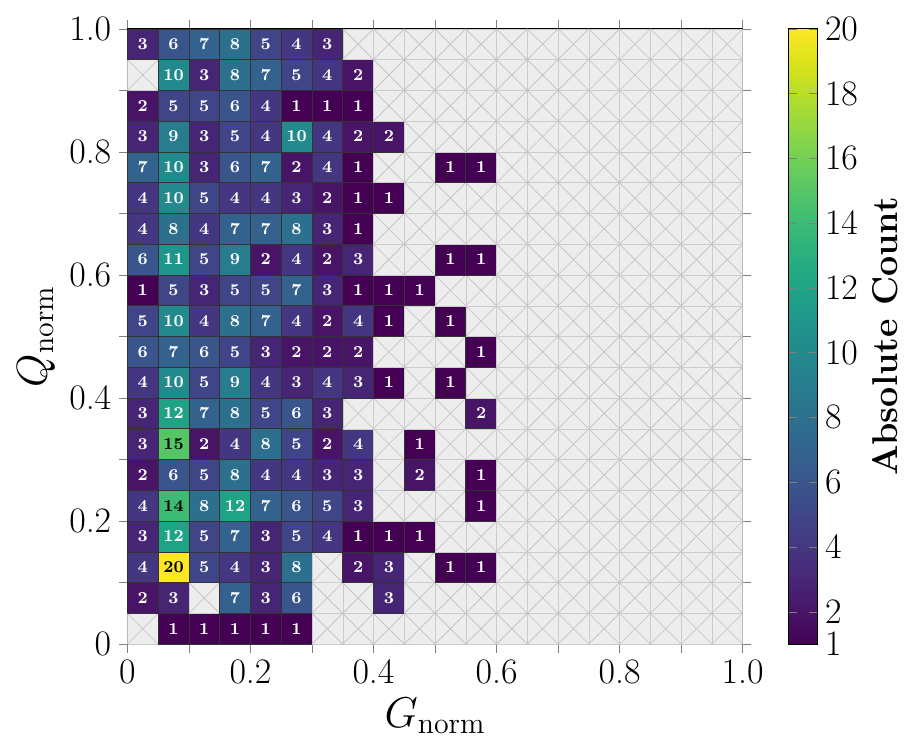}
        \par\smallskip
        \footnotesize\textit{(b) $G_{\text{norm}}$-$Q_{\text{norm}}$ parameter space.}
    \end{minipage}
    \caption{ Data density heatmaps of the generated dataset (in total $770$ simulations) with respect to the parameters $G_{\text{norm}}$-$Q_{\text{norm}}$ and $S_{\text{norm}}$-$Q_{\text{norm}}$.}
    \label{fig:heatmaps}
\end{figure}

To characterize the generated dataset, \cref{fig:heatmaps} shows heatmaps of the sample density in the parameter spaces $S_{\text{norm}}$-$Q_{\text{norm}}$ (a) and $G_{\text{norm}}$-$Q_{\text{norm}}$ (b), 
where lighter colors indicate a higher number of samples per bin.
In the space $S_{\text{norm}}$-$Q_{\text{norm}}$, the data cover a wide range, with the majority of samples falling within the range $0.4 \leq S_{\text{norm}} \leq 0.9$ and $0.1 \leq Q_{\text{norm}} \leq 1.0$. 
Notably, $S_{\text{norm}}$ captures the maximum severity of stenosis across all segments of each geometry, which explains the lower sample density in the left part of the parameter space. 
In the $G_{\text{norm}}$-$Q_{\text{norm}}$ parameter space the data is mostly concentrated in a plane bounded by $G_{\text{norm}} \leq 0.35$, indicating that steeper constrictions, that is, highly severe stenoses across short segments, occur less frequently or not at all in the generated dataset.
The highest sample density occurs in the lower left region of this parameter space.

\subsubsection{Grid Sampling}
\label{postprocessing}

 For CNNs the data must be represented as structured pixel grids. Therefore, the simulation data stored on unstructured grids must first be sampled on pixel grids. An in-house post-processing pipeline based on the Python package \texttt{PyVTK} was implemented for this sampling process \cite{vtkBook}. 
 The pipeline enables the efficient use of \texttt{VTK} operations for sampling the data and computing additional features such as the SDF. The inputs for this pipeline are the XDMF/HDF5 files containing the stored velocity solutions and boundary labels to distinguish between inlet, outlet, and walls of the geometries. The pipeline generates two folders for each simulation: an input folder with the sampled SDF field and an output folder containing the sampled velocities, with the data stored in the CSV file format.

\begin{figure}[ht!]
    \resizebox{\textwidth}{!}{ 
\begin{tikzpicture}
    \usetikzlibrary{decorations.pathreplacing}
            
    \node[anchor=south west,inner sep=0] (image) at (0,0) 
        {\includegraphics[width=\textwidth]{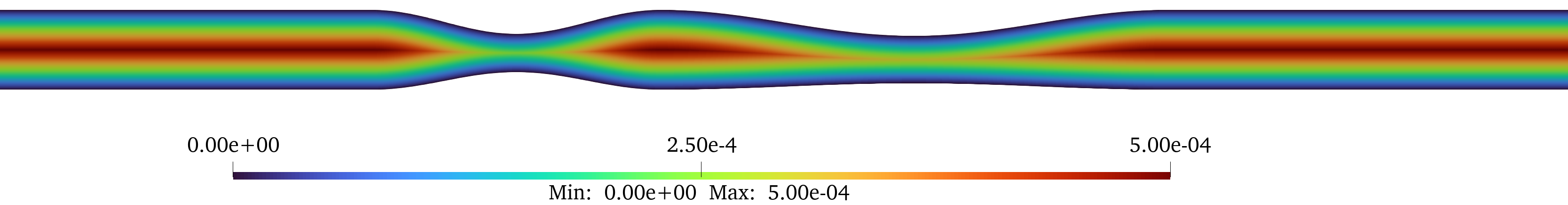}};
       
    \node[above] at (5.9, 0.52) {\scalebox{.5}{\large SDF $[\text{m}]$}}; 

\end{tikzpicture}
}
\vspace{-18pt}
\caption{Computed SDF for an example arterial geometry based on defined boundary labels.}
\label{fig:sdf}
\end{figure}

\paragraph{\textbf{Computation of the SDF}}
\label{sdf}
The main function for computing the SDF is the filter \texttt{vtkImplicitPolyDataDistance} from the \texttt{PyVTK} package. Using the boundary labels of the geometry as input, the filter is applied such that the resulting SDF satisfies equation \cref{eq:sdf}; see \cref{fig:sdf}.

\begin{figure}[ht!]
    \resizebox{\textwidth}{!}{ 
\begin{tikzpicture}
    \usetikzlibrary{decorations.pathreplacing}

    \node[anchor=south west,inner sep=0] (image) at (0,0) 
        {\includegraphics[width=\textwidth,trim={0 100 0 0},clip]{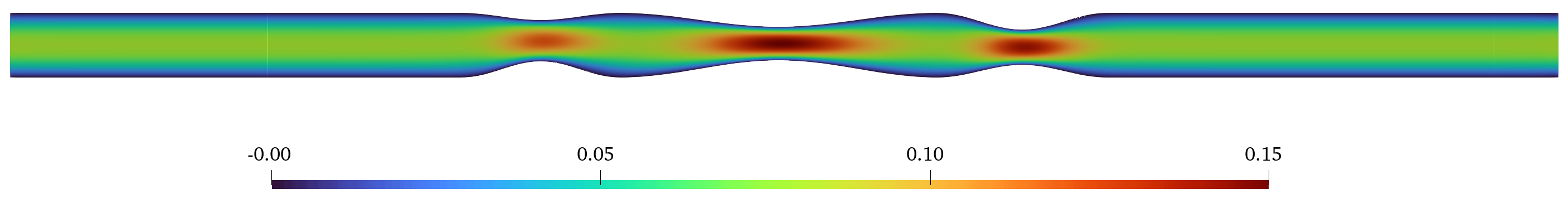}};
    
    \def\xStart{2.24}  
    \def\yMin{1.07}   
    \def\yMax{1.6}   
    \def\xEnd{12.4} 
    \def\yMin{0.58}   
    \def\yMax{1.11}  
    \def\gridXCells{216}  
    \def\gridYCells{12}

    \draw[black, opacity=0.8, line width=1.0pt] 
        (\xStart, \yMin) rectangle (\xEnd, \yMax);
        
    \pgfmathsetmacro\cellWidth{(\xEnd-\xStart)/\gridXCells}
    \pgfmathsetmacro\cellHeight{(\yMax-\yMin)/\gridYCells}
    
    \foreach \x in {0,...,\numexpr\gridXCells-1\relax} {
        \foreach \y in {0,...,\numexpr\gridYCells-1\relax} {
            \pgfmathsetmacro\cellXLeft{\xStart + \x*\cellWidth}
            \pgfmathsetmacro\cellYBottom{\yMin + \y*\cellHeight}
            \pgfmathsetmacro\cellXRight{\cellXLeft + \cellWidth}
            \pgfmathsetmacro\cellYTop{\cellYBottom + \cellHeight}
            
            \draw[black, opacity=0.2, line width=0.5pt] 
                (\cellXLeft, \cellYBottom) rectangle (\cellXRight, \cellYTop);
        }
    }
    
    \foreach \i in {1,...,9} {
        \pgfmathsetmacro\xPos{2.24 + \i * 1.129} 
        \draw[blue, opacity=0.8, line width=1.2pt] (\xPos, \yMin) -- (\xPos, \yMax);
    }

    \foreach \i in {1,...,6} {
        \pgfmathsetmacro\xPos{2.24+0.3 + \i * 1.129}
        \draw[purple,dashed, opacity=0.8, line width=1.2pt] (\xPos, \yMin) -- (\xPos, \yMax);
    }

    \draw[decorate, decoration={brace, amplitude=5pt}, blue, line width=0.6pt] 
        (\xStart+1.129, 0.52) --  (\xStart, 0.52) ;

    \draw[decorate, decoration={brace, amplitude=5pt}, black, line width=0.6pt] 
        (\xStart, \yMax+0.08) -- (\xEnd, \yMax+0.08);

    \pgfmathsetmacro\bracketMid{(\xStart + \xEnd) / 2}
    \node[black] at (\bracketMid, \yMax+0.5) {\scalebox{.8}{{ Global pixel grid}}};

    \pgfmathsetmacro\bracketMid{(\xStart + 1.129/2)}
    \node[black] at (\bracketMid, 0.2) {\scalebox{.8}{{ Subdomain}}};

    \draw[decorate, decoration={brace, amplitude=5pt}, purple, line width=0.6pt] 
        (\xStart+ 0.3 + 1.129 + 1.129 + 1.129, 0.52) --  (\xStart+ 0.3  +1.129 + 1.129, 0.52) ;

    \pgfmathsetmacro\bracketMid{(\xStart + 0.3 + 1.129+ 1.129/2  +1.129  )}
    \node[black] at (\bracketMid, 0.2) {\scalebox{.8}{{ Shifted subdomain}}};

\draw[green!70!black, line width=2pt] ([xshift=2.25cm, yshift=-0.1cm]image.north west) -- ([xshift=2.25cm, yshift=-0.63cm]image.north west) node[pos=0, above, yshift=0.1cm] {$x=x_s$};   

\draw[red, line width=2pt] ([xshift=12.4cm, yshift=-0.1cm]image.north west) -- ([xshift=12.4cm, yshift=-0.63cm]image.north west) node[pos=0, above, yshift=0.1cm] {$x=x_e$};

\end{tikzpicture}
}
\vspace{-20pt}
\caption{Velocity and SDF values are sampled on a global pixel grid. 
The grid starts at $x_s$ and ends at $x_e$ to have developed flow profiles in these cross-sections.
 The resolution is $2305 \times 128$ pixels. Subdomains with $256 \times 128$ pixels are extracted from the global pixel grid for training the local CNN-based USDS. First, $9$ non-overlapping subdomains (solid lines) are extracted. Additionally, $4$ sets of $5$ subdomains are extracted from the stenotic region, with each set being shifted by $50$ pixels in the $x$-direction. The first set of $5$ shifted subdomains is highlighted by the dashed lines.}
\label{fig:global_image_grid}
\end{figure}

\paragraph{\textbf{Sampling of Simulation Data on Structured Pixel Grids}}
 The velocity and SDF are sampled on a global pixel grid measuring $2305 \times 128$ pixels.
 To ensure that fully developed flow profiles are obtained in the straight channel cross-sections at the inlet and outlet, we extract the global pixel grid starting at $x_s = 0.4$ [cm] and ending at $x_e = 2.2$ [cm]; see \cref{fig:global_image_grid}.
  In total, we then obtain $770$ global pixel grids. 

\paragraph{\textbf{Extraction of Subdomains for Training}}
 For training the USDS, a diverse dataset of subdomains with $256 \times 128$ pixels must be extracted, capturing different artery segments with varying flow rates.
 Given $d_{\text{artery}}=1$ [mm], each subdomain captures an artery segment with dimensions of $2$ [mm] $\times$ $1$ [mm]. 
First, $9$ non-overlapping subdomains are extracted from each global pixel grid; see \cref{fig:global_image_grid}. The dataset is then further augmented by extracting $5\cdot4=20$ additional subdomains from the stenotic region. Each set of $5$ subdomains is shifted by $50$ pixels in the $x$-direction in order to capture different segments. Adding these additional subdomains enriches the training data and should prevent the model from learning only the non-overlapping decomposition pattern. In total, $(9+5\cdot4)\cdot770 = 22.330$ training images are obtained.

\subsection{Training Details}
\label{training}
The settings used for the optimization of the CNN network model are described here. The SDF data is scaled to the range $[0,1]$ by applying max normalization, where all values are divided by $\text{SDF}_{\max}= 5.0\cdot 10^{-4}$.
For the velocity data, the best results were currently achieved without normalization. \cref{tab:data_statistics} shows the statistical properties of the data.
Notably, the scale difference between the $x$- and $y$-component  is apparent, as the $x$-component is dominating the flow characteristics for the channel-like flows.
 \begin{table}[htbp]
\centering
\small
\begin{tabular}{l|ccc}
\toprule
& SDF & $x$-component \textbf{$v_x$} &  $y$-component  \textbf{$v_y$} \\
\midrule
\textbf{Mean} & $3.53 \cdot 10^{-1}$ & $1.97 \cdot 10^{-1}$ & $1.66 \cdot 10^{-6}$ \\
\textbf{Standard deviation} & $2.89 \cdot 10^{-1}$ & $2.00 \cdot 10^{-1}$ & $1.07 \cdot 10^{-2}$ \\
\bottomrule
\end{tabular}
\vspace{5pt}
\caption{Mean and standard deviation of whole dataset with normalized SDF}
\label{tab:data_statistics}
\end{table}

As loss function $\mathcal{L}$, the mean squared error (MSE) is used with an additional penalization term for the scaling factors $C$ in the constraint layer, that is,
\begin{equation}
    \label{eq:loss}
    \mathcal{L} = \frac{1}{|D|} \sum_{n \in D} \left( \frac{1}{|I_n|} \sum_{p \in I_n} \| \bm{v}^{\text{FEM}}_p - \bm{v}^{\text{CNN}}_p \|_2^2  + \lambda_C  \frac{1}{W} \sum_{i=0}^{W-1} \left( C_n(i) - 1 \right)^2    \right),
\end{equation}
where $D$ represents a set of images, $I_n$ is the set of pixels in image $n$, $\bm{v}^{\text{FEM}}_p $ is the FEM velocity at pixel $p$, $\bm{v}^{\text{CNN}}_p$ is the predicted velocity at pixel $p$,
$\lambda_C$ is a weighting factor and $C_n(i)$ denotes the calculated scaling factor at pixel column $i$ for image $n$.
For the purely data-driven CNN, $\lambda_C = 0$ as the constraint layer is not present.
For the flow-rate-conserving CNN,  training without the penalization term ($\lambda_C = 0$) leads to low loss values, however, 
 the resulting models generate scaling factors with an average of  $\Bar{C} \approx 0.2$.
This indicates that the network learns to overestimate the velocity magnitudes before scaling, which are then downscaled by the constraint layer.
Setting $\lambda_C = 10^{-3}$ results in models with low loss values and whose scaling factors remain close to one,
such that the CNN learns quantitatively accurate velocity fields also before scaling,
and the constraint layer makes only minor corrections as intended. Larger values $\lambda_C  > 10^{-2}$ prevented the training from converging.
We therefore use $\lambda_C = 10^{-3}$ for the training of all flow-rate-conserving models.

As the velocity components have different scales, the MSE loss is dominated by the $x$-component. 
To improve the prediction accuracy for $v_y$, a weighted MSE loss was investigated.
While it improved the prediction accuracy for the $y$-component for subdomain prediction,  it led to
a slight increase in the $v_x$ error. Using these models as USDS in the CNN-Schwarz-Flow approach resulted in overall higher global prediction errors, as $v_x$ is the dominant flow component.
We therefore use the loss \cref{eq:loss} for this study. Further investigations regarding the loss function design and normalization strategies are planned for future work.

The dataset is split into training, validation, and test sets with a ratio of $80\%$, $10\%$, and $10\%$, respectively. 
Batch sizes of $64$ for training and $32$ for validation and testing are used.
We employ the Adam (Adapative Momentum Method) optimizer from the Python package \texttt{PyTorch}
with an initial learning rate of $\alpha_0 = 2.5 \cdot 10^{-4}$  \cite{pytorch, kingma2017adammethodstochasticoptimization}.
The \texttt{ReduceLROnPlateau} scheduler is applied to halve the learning rate if the validation loss does not improve by more than $10^{-7}$ over $20$ consecutive epochs, down to a minimum of $\alpha_{\min} = 10^{-5}$.
Furthermore, an early stopping mechanism is used, which is a form of regularization technique to prevent overfitting and terminate the training if the validation loss no longer significantly improves \cite{prechelt2002early}. 
In this context, we set the patience to $50$ epochs, the minimum improvement value to $10^{-7}$, and the maximum number of training epochs to $600$.
The CNNs and the entire training framework are implemented with \texttt{PyTorch} functions, using, for example, the default settings for weight initialization. 
By leveraging the \texttt{Horovod} Python package,  we implemented a data-parallel approach to distribute the workload across multiple GPUs, with each GPU processing a portion of data and synchronizing gradients during backpropagation \cite{horovod}. 
The training of the CNN surrogate models were then performed on the Wisteria/BDEC-01 Aquarius supercomputer, employing $5$ GPUs, where $600$ epochs took approximately $2$ hours to train.
 For reproducibility and parameter tracking, we use the \texttt{Hydra} Python framework  \cite{yadan2019hydra}.

\section{Performance Analysis}
\label{sec:results}

First, we focus on the USDS and examine the training and generalization performance for a single subdomain prediction, focusing particulary on how the width of the input boundary $\xi
$ and the CNN architecture affect model convergence and inference accuracy. Second, applying the 
 CNN-Schwarz-Flow approach using the trained USDS models, we assess and compare the accuracies for the global prediction. 
Finally, we investigate the transferability of the best performing model to arterial domains of increasing length.

\subsection{Subdomain Prediction: Training of CNN Surrogate Models}
\label{single_pred_performance}

 We train both the data-driven model and the flow-rate-conserving model for three input boundary
widths $\xi=\{1,10,20\}$, resulting in a total of six different CNNs.
We restrict the analysis to these values,  as larger $\xi$ increases the number of pixel columns overwritten by the boundary enforcement and, 
due to the overlap condition $\delta \geq 2\xi$, requires more subdomains to cover the full domain.
Due to the inherent stochasticity of our GPU-accelerated training, the loss curves for identical hyperparameters can vary between individual runs, but their overall trend remains consistent. 
We performed multiple training runs and selected the final models based on the minimum achieved losses.
An overview and comparison of loss values achieved by the selected models are shown in  \cref{fig:combined_single_pred}.
  \begin{figure}[ht!]
    \centering
        \begin{minipage}{0.49\textwidth}
        \centering
        \resizebox{\textwidth}{!}{
%
%
\definecolor{mycolor1}{rgb}{0.40000,0.80000,0.40000}%
\definecolor{mycolor2}{rgb}{0.80000,0.40000,0.80000}%
\begin{tikzpicture}

\begin{axis}[%
width=5.5in,
height=3.0in,
at={(1.156in,0.965in)},
scale only axis,
xmin=0.5,
xmax=9.5,
xtick={1,2,3,4,5,6,7,8,9},
xticklabels={{Train.},{Val.},{Test},{Train.},{Val.},{Test},{Train.},{Val.},{Test}},
yticklabel style={font=\LARGE},
xticklabel style={rotate=50, font=\LARGE},
xtick style={font=\LARGE},
ytick style={font=\LARGE},
xlabel={\textbf{Dataset}},
ymode=log,
ymin=1e-07,
ymax=1e-05,
yminorticks=true,
xlabel style={font=\LARGE\color{white!15!black}, yshift=-30pt},
ylabel style={font=\LARGE\color{white!15!black}, yshift=20pt},
ylabel={\textbf{MSE loss}},
axis background/.style={fill=white},
title style={font=\LARGE},
title={\textbf{Flow-rate-conserving CNN}},
axis x line*=bottom,
axis y line*=left,
xmajorgrids,
ymajorgrids,
yminorgrids,
legend style={at={(0.03,0.19)}, anchor=west, legend cell align=left, align=left, draw=white!15!black, font=\LARGE}
]
\addplot [color=white!10!orange, line width=3pt, mark=diamond*, mark options={line width=1pt}, mark size=4pt]
  table[row sep=crcr]{%
1	1.32e-06\\
2	2.22e-06\\
3	2.06e-06\\
};
\addlegendentry{$\xi=1$}

\addplot [color=white!10!orange, line width=3pt, forget plot]
  table[row sep=crcr]{%
4	8e-07\\
5	1.55e-06\\
6	1.44e-06\\
};
\addplot [color=white!10!orange, line width=3pt, forget plot]
  table[row sep=crcr]{%
7	5.2e-07\\
8	6.73e-07\\
9	6.23e-07\\
};
\addplot[only marks, mark=diamond*, mark options={}, mark size=4pt, color=white!10!orange, fill=white!10!orange, forget plot] table[row sep=crcr]{%
x	y\\
1	1.32e-06\\
2	2.22e-06\\
3	2.06e-06\\
4	8e-07\\
5	1.55e-06\\
6	1.44e-06\\
7	5.2e-07\\
8	6.73e-07\\
9	6.23e-07\\
};

\addplot [color=white!10!mycolor1, line width=3pt, mark=x, mark options={line width=2pt}, mark size=6pt]
  table[row sep=crcr]{%
1	1.15e-06\\
2	2.01e-06\\
3	1.83e-06\\
};
\addlegendentry{$\xi=10$}

\addplot [color=white!10!mycolor1, line width=3pt, forget plot]
  table[row sep=crcr]{%
4	7.27e-07\\
5	1.43e-06\\
6	1.31e-06\\
};
\addplot [color=white!10!mycolor1, line width=3pt, forget plot]
  table[row sep=crcr]{%
7	4.26e-07\\
8	5.81e-07\\
9	5.24e-07\\
};
\addplot[only marks, mark=x, mark options={line width=2pt}, mark size=6pt, color=white!10!mycolor1, fill=white!10!mycolor1, forget plot] table[row sep=crcr]{%
x	y\\
1	1.15e-06\\
2	2.01e-06\\
3	1.83e-06\\
4	7.27e-07\\
5	1.43e-06\\
6	1.31e-06\\
7	4.26e-07\\
8	5.81e-07\\
9	5.24e-07\\
};

\addplot [color=white!10!blue, line width=3pt, mark=square*, mark options={line width=1pt}, mark size=3pt]
  table[row sep=crcr]{%
1	8.41818746266654e-07\\
2	1.61891182415275e-06\\
3	1.50361393025378e-06\\
};
\addlegendentry{$\xi=20$}

\addplot [color=white!10!blue, line width=3pt, forget plot]
  table[row sep=crcr]{%
4	5.72986948554899e-07\\
5	1.2222966461195e-06\\
6	1.14351314550731e-06\\
};
\addplot [color=white!10!blue, line width=3pt, forget plot]
  table[row sep=crcr]{%
7	2.68831769290045e-07\\
8	3.96615035924697e-07\\
9	3.60100870011593e-07\\
};
\addplot[only marks, mark=square*, mark options={}, mark size=3pt, color=white!10!blue, fill=white!10!blue, forget plot] table[row sep=crcr]{%
x	y\\
1	8.41818746266654e-07\\
2	1.61891182415275e-06\\
3	1.50361393025378e-06\\
4	5.72986948554899e-07\\
5	1.2222966461195e-06\\
6	1.14351314550731e-06\\
7	2.68831769290045e-07\\
8	3.96615035924697e-07\\
9	3.60100870011593e-07\\
};

\addplot [color=gray, dashed, forget plot]
  table[row sep=crcr]{%
3.5	1e-07\\
3.5	1e-05\\
};
\addplot [color=gray, dashed, forget plot]
  table[row sep=crcr]{%
6.5	1e-07\\
6.5	1e-05\\
};
\node[draw, fill=white, opacity=0.8, text opacity=1, rounded corners, align=left, font=\LARGE] 
  at (axis cs:2,8e-6) [anchor=center] 
  {
    Velocity $\bm{v}$
  };
  
  \node[draw, fill=white, opacity=0.8, text opacity=1, rounded corners, align=left, font=\LARGE] 
  at (axis cs:5,8e-6) [anchor=center] 
  {
    $x$-component $v_x$
  };
  
    \node[draw, fill=white, opacity=0.8, text opacity=1, rounded corners, align=left, font=\LARGE] 
  at (axis cs:8,8e-6) [anchor=center] 
  {
    $y$-component $v_y$
  };
\end{axis}
\end{tikzpicture}%
}
        \label{fig:figure2_hc_singlepred}
    \end{minipage}
    \hfill
        \begin{minipage}{0.49\textwidth}
        \centering
        \resizebox{\textwidth}{!}{
%
%
\definecolor{mycolor1}{rgb}{0.40000,0.80000,0.40000}%
\definecolor{mycolor2}{rgb}{0.80000,0.40000,0.80000}%
\begin{tikzpicture}

\begin{axis}[%
width=5.5in,
height=3.0in,
at={(1.156in,0.965in)},
scale only axis,
xmin=0.5,
xmax=9.5,
xtick={1,2,3,4,5,6,7,8,9},
xticklabels={{Train.},{Val.},{Test},{Train.},{Val.},{Test},{Train.},{Val.},{Test}},
yticklabel style={font=\LARGE},
xticklabel style={rotate=50, font=\LARGE},
xtick style={font=\LARGE},
ytick style={font=\LARGE},
xlabel={\textbf{Dataset}},
ymode=log,
ymin=1e-07,
ymax=1e-05,
yminorticks=true,
xlabel style={font=\LARGE\color{white!15!black}, yshift=-30pt},
ylabel style={font=\LARGE\color{white!15!black}, yshift=20pt},
ylabel={\textbf{MSE loss}},
axis background/.style={fill=white},
title style={font=\LARGE},
title={\textbf{Data-driven CNN}},
axis x line*=bottom,
axis y line*=left,
xmajorgrids,
ymajorgrids,
yminorgrids,
legend style={at={(0.03,0.19)}, anchor=west, legend cell align=left, align=left, draw=white!15!black, font=\LARGE}
]
\addplot [color=white!40!orange, dashed, line width=3pt, mark=diamond*, mark options={solid, line width=1pt}, mark size=4pt]
  table[row sep=crcr]{%
1	1.62e-06\\
2	2.72e-06\\
3	2.75e-06\\
};
\addlegendentry{$\xi=1$ }

\addplot [color=white!40!orange, dashed, line width=3pt, forget plot]
  table[row sep=crcr]{%
4	1.13e-06\\
5	2.07e-06\\
6	2.13e-06\\
};
\addplot [color=white!40!orange, dashed, line width=3pt, forget plot]
  table[row sep=crcr]{%
7	4.85e-07\\
8	6.49e-07\\
9	6.23e-07\\
};
\addplot[only marks, mark=diamond*, mark options={}, mark size=4pt, color=white!40!orange, fill=white!40!orange, forget plot] table[row sep=crcr]{%
x	y\\
1	1.62e-06\\
2	2.72e-06\\
3	2.75e-06\\
4	1.13e-06\\
5	2.07e-06\\
6	2.13e-06\\
7	4.85e-07\\
8	6.49e-07\\
9	6.23e-07\\
};
\addplot [color=mycolor1, dashed,  line width=3pt, mark=x, mark options={solid, line width=2pt}, mark size=6pt]
  table[row sep=crcr]{%
1	1.23e-06\\
2	2.32e-06\\
3	2.12e-06\\
};
\addlegendentry{$\xi=10$ }

\addplot [color=mycolor1, dashed, line width=3pt, forget plot]
  table[row sep=crcr]{%
4	9.06e-07\\
5	1.86e-06\\
6	1.74e-06\\
};
\addplot [color=mycolor1, dashed, line width=3pt, forget plot]
  table[row sep=crcr]{%
7	3.69e-07\\
8	5.19e-07\\
9	4.8e-07\\
};
\addplot[only marks,mark=x, mark options={line width=2pt}, mark size=6pt, color=mycolor1, fill=mycolor1, forget plot] table[row sep=crcr]{%
x	y\\
1	1.23e-06\\
2	2.32e-06\\
3	2.12e-06\\
4	9.06e-07\\
5	1.86e-06\\
6	1.74e-06\\
7	3.69e-07\\
8	5.19e-07\\
9	4.8e-07\\
};
\addplot [color=white!40!blue, dashed,  line width=3pt, mark=square*, mark options={solid, line width=1pt}, mark size=3pt]
  table[row sep=crcr]{%
1	1.2e-06\\
2	2.25e-06\\
3	1.98e-06\\
};
\addlegendentry{$\xi=20$}

\addplot [color=white!40!blue, dashed, line width=3pt, forget plot]
  table[row sep=crcr]{%
4	8.74e-07\\
5	1.82e-06\\
6	1.65e-06\\
};
\addplot [color=white!40!blue, dashed, line width=3pt, forget plot]
  table[row sep=crcr]{%
7	3.6e-07\\
8	5.03e-07\\
9	4.69e-07\\
};
\addplot[only marks, mark=square*, mark options={}, mark size=3pt, color=white!40!blue, fill=white!40!blue, forget plot] table[row sep=crcr]{%
x	y\\
1	1.2e-06\\
2	2.25e-06\\
3	1.98e-06\\
4	8.74e-07\\
5	1.82e-06\\
6	1.65e-06\\
7	3.6e-07\\
8	5.03e-07\\
9	4.69e-07\\
};
\addplot [color=gray, dashed, forget plot]
  table[row sep=crcr]{%
3.5	1e-07\\
3.5	1e-05\\
};
\addplot [color=gray, dashed, forget plot]
  table[row sep=crcr]{%
6.5	1e-07\\
6.5	1e-05\\
};

\node[draw, fill=white, opacity=0.8, text opacity=1, rounded corners, align=left, font=\LARGE] 
  at (axis cs:2,8e-6) [anchor=center] 
  {
    Velocity $\bm{v}$
  };
  
  \node[draw, fill=white, opacity=0.8, text opacity=1, rounded corners, align=left, font=\LARGE] 
  at (axis cs:5,8e-6) [anchor=center] 
  {
    $x$-component $v_x$
  };
  
    \node[draw, fill=white, opacity=0.8, text opacity=1, rounded corners, align=left, font=\LARGE] 
  at (axis cs:8,8e-6) [anchor=center] 
  {
    $y$-component $v_y$
  };
\end{axis}
\end{tikzpicture}%
}
        \label{fig:figure1_dd_singlepred}
    \end{minipage}
    \vspace{-17pt}
    \caption{ Comparison of achieved MSE losses for the trained data-driven CNN (right) and flow-rate-conserving CNN architecture (left) using $\xi=\{1,10,20\}$.
    The results are divided into MSE losses for the velocity vector $\bm{v}$ and for the individual $x$- and $y$-velocitiy components. The loss values are computed for training (Train.), validation (Val.) and test (Test) datasets. }
    \label{fig:combined_single_pred}
\end{figure}
  The plot shows for training, validation, and test datasets, the MSE loss for the velocity vector $\boldsymbol{v}$ and its components $v_x$ and $v_y$ at the determined final epoch for each model.
   Since the $y$-component values are at least an order of magnitude smaller than the $x$-component values, the MSE losses for $v_y$ are respectively smaller.
   The flow-rate-conserving models achieve slightly lower loss values than the data-driven models; 
   for example, for $\xi=1$ the validation MSE for the velocity vector $\bm{v}$ is $2.2 \cdot 10^{-6}$ for the constrained model and $2.7 \cdot 10^{-6}$ for the data-driven model. 
   Increasing $\xi$ reduces the losses only marginally for both models.
   Note that larger $\xi$ values also increase the number of boundary columns that are overwritten with the target FEM values, which are accounted for in the MSE loss.
   Based on these results, we note that, given the present dataset, architecture, and training configurations, $\xi=1$ is sufficient to achieve low loss values for both architectures.

\subsection{Global Prediction: Application of the CNN-Schwarz-Flow Approach}
\label{global_pred_performance}

The six trained surrogate models from the previous section are now each used as inexact subdomain solvers within the 
CNN-Schwarz-Flow approach; see \cref{sec:global_image_prediction}.
 The entire configuration details are listed in \cref{tab:dd_combined_params}.
\begin{table}[htbp]
\centering
\small  
\renewcommand{\arraystretch}{1.2}  
\resizebox{\textwidth}{!}{
\begin{tabular}{c c c c c c}
\toprule[1.0pt]
\multicolumn{5}{c}{\textbf{CNN-Schwarz-Flow Algorithm Parameters} } \\
\midrule
\textbf{Overlap width} $\delta$ & \textbf{Initialization}  & \textbf{Stopping criterion} & \textbf{Tolerance} $\epsilon$ & \textbf{Max. iteration} \\
 $2\xi$ & Parabolic  & $\text{AbsErr}_{\max}$ \cref{eq:max:abs_error} & $1.0 \cdot 10^{-4}$ & $200$ \\
\midrule[1.0pt]
\multicolumn{5}{c}{\textbf{Domain Decomposition Configuration}} \\
\midrule
\textbf{$\xi$/$\delta$} & $\xi=1$/$\delta=2$ & $\xi=10$/$\delta=20$ & $\xi=20$/$\delta=40$ &  \\
\textbf{$\#$Overlapping subdomains} & $9$ & $9$& $10$  & \\
\textbf{Covered pixel range (Width)} & \( [0, 2288] \) & \( [0, 2144] \) & \( [0, 2200] \) & \\
\bottomrule[1.0pt]
\end{tabular}
}
\caption{Prescribed parameters in the CNN-Schwarz-Flow algorithm.}
\label{tab:dd_combined_params}
\end{table}
We first test the algorithm for all $770$ simulations on the corresponding global pixel grids of size $2305 \times 128$ pixels.
Each global grid is decomposed into overlapping subdomains with a respective overlap width of $\delta=2\xi$.  The overlap $\delta$ determines both the number of overlapping subdomains and the total pixel area covered by the union of the subdomains. 
The extended outlet area ensures that fully developed flow profiles are maintained in the outlet cross-section, making the slight deviations in pixel coverage not relevant for further analysis. 
For both the flow-rate-conserving and data-driven CNNs the velocity $x$-component for the interior subdomains is initialized with a parabolic velocity profile; see \cref{sec:global_image_prediction}.
 We use an empirically determined value for the tolerance $\epsilon$. In general, the choice of the stopping criterion and the tolerance remains an open research question and depends, for example, on the specific problem scales.

As an error metric, we compute the global relative error (GRE), that is, 
 \begin{equation}
\label{eq:gre}
\text{GRE} = \frac{\left\| \bm{v} - \hat{\bm{v}} \right\|_2}{\left\| \bm{v} \right\|_2 + 10^{-4}},
\end{equation} 
where $\bm{v}$ is the reference velocity solution computed from the FEM simulation and $\hat{\bm{v}}$ is the predicted velocity field using the CNN-Schwarz-Flow approach.
Notably, only pixels within the geometry are considered for computing the GRE by applying the mask \cref{eq:mask}.
Furthermore, we define a lower bound for the error,  denoted by GRE$^*$, by fixing the target FEM velocities as input velocities 
for each subdomain and computing the GRE after a single red-black iteration.
Thus, GRE$^*$ quantifies the error introduced by the CNN being an inexact subdomain solver, independent of the Schwarz iterations.

In  \cref{fig:combined_multipred}, two histograms comparing the achieved GRE values for the flow-rate-conserving and data-driven models applied to the entire dataset are shown.
\definecolor{mycolor2}{rgb}{0.40000,0.80000,0.40000}
\pgfplotstableread[row sep=\\,col sep=&]{
    category & xi1 & xi10 & xi20 \\
    {$\leq 1\%$} & 644  & 560 & 631 \\
    {1-5\%}      & 122  & 197 & 135 \\
    {5-10\%}     & 4     & 10  & 4   \\
    {10-20\%}    & 0      & 1   & 0   \\
    {$>$20\%}    & 0      & 1   & 0   \\
    {diverged}   & 0      & 1   & 0   \\
}\mydatahc
\pgfplotstableread[row sep=\\,col sep=&]{
    category & xi1  & xi10 & xi20 \\
    {$\leq 1\%$} & 307  & 139 & 72  \\
    {1-5\%}      & 427  & 507 & 414 \\
    {5-10\%}     & 21   & 78  & 132 \\
    {10-20\%}    & 12  & 23  & 86  \\
    {$>$20\%}    & 3     & 23  & 66  \\
    {diverged}   & 0     & 0   & 0   \\
}\mydatadatadriven
\begin{figure}[h!]
    \vspace{-10pt}
    \begin{minipage}{\textwidth}
        \centering
        \resizebox{\textwidth}{!}{\begin{tikzpicture}
    \begin{axis}[
            ybar,
            bar width=.3cm,
            width=\textwidth,
            height=.27\textwidth,
            legend style={at={(0.65,1)},
            anchor=north,legend columns=-1, font=\scriptsize},
            symbolic x coords={{$\leq 1\%$},{1-5\%},{5-10\%},{10-20\%},{$>$20\%},{diverged}},
            xtick=data,
            ymin=0.0,ymax=770,
            ylabel={\textbf{\# Problems}},
            xlabel={\textbf{Global relative error (GRE)}},
            x label style={yshift=0.0cm, font=\scriptsize},
            y label style={font=\scriptsize, yshift=-0.35cm},
            title style={font=\footnotesize},
            title={\textbf{USDS: Flow-rate-conserving CNN}},
            x tick label style={xshift=0.5cm, yshift=-0.2cm, anchor=east, font=\scriptsize},
            ytick={0,200,400,600,770},
            y tick label style={font=\scriptsize},
        ]
        \addplot[fill=white!10!orange!70, draw=white!10!orange!70, draw=none ] table[x=category,y=xi1]{\mydatahc};
        \addplot[fill=white!10!mycolor2,   draw=white!10!mycolor2]   table[x=category,y=xi10]{\mydatahc};
        \addplot[fill=white!10!blue!70,    draw=white!10!blue!70, draw=none]    table[x=category,y=xi20]{\mydatahc};

        \draw[gray, dashed, xshift=-0.9cm] ({axis cs:{1-5\%},0}   |- {rel axis cs:0,0}) -- ({axis cs:{1-5\%},0}   |- {rel axis cs:0,1});
        \draw[gray, dashed, xshift=-0.9cm] ({axis cs:{5-10\%},0}  |- {rel axis cs:0,0}) -- ({axis cs:{5-10\%},0}  |- {rel axis cs:0,1});
        \draw[gray, dashed, xshift=-0.9cm] ({axis cs:{10-20\%},0} |- {rel axis cs:0,0}) -- ({axis cs:{10-20\%},0} |- {rel axis cs:0,1});
        \draw[gray, dashed, xshift=-0.9cm] ({axis cs:{$>$20\%},0} |- {rel axis cs:0,0}) -- ({axis cs:{$>$20\%},0} |- {rel axis cs:0,1});
        \draw[gray, dashed, xshift=-0.9cm] ({axis cs:{diverged},0}|- {rel axis cs:0,0}) -- ({axis cs:{diverged},0}|- {rel axis cs:0,1});

        \node[rotate=90, font=\scriptsize, yshift=11, xshift=-12] at (axis cs:{$\leq 1\%$},644) {83.6\%};
        \node[rotate=90, font=\scriptsize, yshift=11, xshift=12]  at (axis cs:{1-5\%},122)       {15.8\%};
        \node[rotate=90, font=\scriptsize, yshift=11, xshift=12]  at (axis cs:{5-10\%},4)         {0.5\%};
        \node[rotate=90, font=\scriptsize, yshift=11, xshift=12]  at (axis cs:{10-20\%},0)        {0.0\%};
        \node[rotate=90, font=\scriptsize, yshift=11, xshift=12]  at (axis cs:{$>$20\%},0)        {0.0\%};
        \node[rotate=90, font=\scriptsize, yshift=11, xshift=12]  at (axis cs:{diverged},0)       {0.0\%};

        \node[rotate=90, font=\scriptsize, yshift=0, xshift=-12] at (axis cs:{$\leq 1\%$},560) {72.7\%};
        \node[rotate=90, font=\scriptsize, yshift=0, xshift=12]  at (axis cs:{1-5\%},197)       {25.6\%};
        \node[rotate=90, font=\scriptsize, yshift=0, xshift=12]  at (axis cs:{5-10\%},10)        {1.3\%};
        \node[rotate=90, font=\scriptsize, yshift=0, xshift=12]  at (axis cs:{10-20\%},1)        {0.1\%};
        \node[rotate=90, font=\scriptsize, yshift=0, xshift=12]  at (axis cs:{$>$20\%},1)        {0.1\%};
        \node[rotate=90, font=\scriptsize, yshift=0, xshift=12]  at (axis cs:{diverged},1)       {0.1\%};

        \node[rotate=90, font=\scriptsize, color=white, yshift=-11, xshift=-12] at (axis cs:{$\leq 1\%$},631) {81.9\%};
        \node[rotate=90, font=\scriptsize, yshift=-11, xshift=12]  at (axis cs:{1-5\%},135)       {17.5\%};
        \node[rotate=90, font=\scriptsize, yshift=-11, xshift=12]  at (axis cs:{5-10\%},4)         {0.5\%};
        \node[rotate=90, font=\scriptsize, yshift=-11, xshift=12]  at (axis cs:{10-20\%},0)        {0.0\%};
        \node[rotate=90, font=\scriptsize, yshift=-11, xshift=12]  at (axis cs:{$>$20\%},0)        {0.0\%};
        \node[rotate=90, font=\scriptsize, yshift=-11, xshift=12]  at (axis cs:{diverged},0)       {0.0\%};

        \legend{$\xi=1$,  $\xi=10$, $\xi=20$}
    \end{axis}
\end{tikzpicture}}
        \vspace{-13pt}
        \label{fig:figure1_hc_multipred}
    \end{minipage}
    \begin{minipage}{\textwidth}
        \centering
        \resizebox{\textwidth}{!}{\begin{tikzpicture}
    \begin{axis}[
            ybar,
            bar width=.3cm,
            width=\textwidth,
            height=.27\textwidth,
            legend style={at={(0.65,1)},
            anchor=north,legend columns=-1, font=\scriptsize},
            symbolic x coords={{$\leq 1\%$},{1-5\%},{5-10\%},{10-20\%},{$>$20\%},{diverged}},
            xtick=data,
            ymin=0.0,ymax=770,
            ylabel={\textbf{\# Problems}},
            xlabel={\textbf{Global relative error (GRE)}},
            x label style={yshift=0.0cm, font=\scriptsize},
            y label style={font=\scriptsize, yshift=-0.35cm},
            title style={font=\footnotesize},
            title={\textbf{USDS: Data-driven CNN}},
            x tick label style={xshift=0.5cm, yshift=-0.2cm, anchor=east, font=\scriptsize},
            ytick={0,200,400,600,770},
            y tick label style={font=\scriptsize},
        ]
        \addplot[fill=white!40!orange!70, draw=white!40!orange!70,  draw=none] table[x=category,y=xi1]{\mydatadatadriven};
        \addplot[fill=white!40!mycolor2,   draw=white!40!mycolor2,   draw=none] table[x=category,y=xi10]{\mydatadatadriven};
        \addplot[fill=white!40!blue!70,    draw=white!40!blue!70,     draw=none] table[x=category,y=xi20]{\mydatadatadriven};

        \draw[gray, dashed, xshift=-0.9cm] ({axis cs:{1-5\%},0}   |- {rel axis cs:0,0}) -- ({axis cs:{1-5\%},0}   |- {rel axis cs:0,1});
        \draw[gray, dashed, xshift=-0.9cm] ({axis cs:{5-10\%},0}  |- {rel axis cs:0,0}) -- ({axis cs:{5-10\%},0}  |- {rel axis cs:0,1});
        \draw[gray, dashed, xshift=-0.9cm] ({axis cs:{10-20\%},0} |- {rel axis cs:0,0}) -- ({axis cs:{10-20\%},0} |- {rel axis cs:0,1});
        \draw[gray, dashed, xshift=-0.9cm] ({axis cs:{$>$20\%},0} |- {rel axis cs:0,0}) -- ({axis cs:{$>$20\%},0} |- {rel axis cs:0,1});
        \draw[gray, dashed, xshift=-0.9cm] ({axis cs:{diverged},0}|- {rel axis cs:0,0}) -- ({axis cs:{diverged},0}|- {rel axis cs:0,1});

        \node[rotate=90, font=\scriptsize, yshift=11, xshift=12]  at (axis cs:{$\leq 1\%$},307)  {39.9\%};
        \node[rotate=90, font=\scriptsize, yshift=11, xshift=12]  at (axis cs:{1-5\%},427)        {55.5\%};
        \node[rotate=90, font=\scriptsize, yshift=11, xshift=12]  at (axis cs:{5-10\%},21)         {2.7\%};
        \node[rotate=90, font=\scriptsize, yshift=11, xshift=12]  at (axis cs:{10-20\%},12)        {1.6\%};
        \node[rotate=90, font=\scriptsize, yshift=11, xshift=12]  at (axis cs:{$>$20\%},3)         {0.4\%};
        \node[rotate=90, font=\scriptsize, yshift=11, xshift=12]  at (axis cs:{diverged},0)        {0.0\%};

        \node[rotate=90, font=\scriptsize, yshift=0, xshift=12]  at (axis cs:{$\leq 1\%$},139)  {18.1\%};
        \node[rotate=90, font=\scriptsize, yshift=0, xshift=-12] at (axis cs:{1-5\%},507)        {65.8\%};
        \node[rotate=90, font=\scriptsize, yshift=0, xshift=12]  at (axis cs:{5-10\%},78)        {10.1\%};
        \node[rotate=90, font=\scriptsize, yshift=0, xshift=12]  at (axis cs:{10-20\%},23)        {3.0\%};
        \node[rotate=90, font=\scriptsize, yshift=0, xshift=12]  at (axis cs:{$>$20\%},23)        {3.0\%};
        \node[rotate=90, font=\scriptsize, yshift=0, xshift=12]  at (axis cs:{diverged},0)        {0.0\%};

        \node[rotate=90, font=\scriptsize, yshift=-11, xshift=12] at (axis cs:{$\leq 1\%$},72)   {9.4\%};
        \node[rotate=90, font=\scriptsize, yshift=-11, xshift=12] at (axis cs:{1-5\%},414)        {53.8\%};
        \node[rotate=90, font=\scriptsize, yshift=-11, xshift=12] at (axis cs:{5-10\%},132)       {17.1\%};
        \node[rotate=90, font=\scriptsize, yshift=-11, xshift=12] at (axis cs:{10-20\%},86)       {11.2\%};
        \node[rotate=90, font=\scriptsize, yshift=-11, xshift=12] at (axis cs:{$>$20\%},66)        {8.6\%};
        \node[rotate=90, font=\scriptsize, yshift=-11, xshift=12] at (axis cs:{diverged},0)        {0.0\%};

        \legend{$\xi=1$, $\xi=10$, $\xi=20$}
    \end{axis}
\end{tikzpicture}}
        \label{fig:figure2_dd_multipred}
    \end{minipage}
    \vspace{-19pt}
    \caption{Histograms of the achieved GRE values for the CNN-Schwarz-Flow algorithm applied to the
             entire dataset ($770$ simulations), using the flow-rate-conserving (upper) and
             data-driven (lower) CNN models from \cref{single_pred_performance}. The x-axis is divided into six categories based on
the GRE \cref{eq:gre}. The y-axis shows the total number of global predictions that fall into each category
using the respective model. On each bar the relative percentage is also displayed. }
    \label{fig:combined_multipred}
\end{figure}
 
For the flow-rate-conserving CNN model, the best results are achieved for $\xi=1$, where $83.6$\% of the predictions have a GRE $\leq 1\%$ and only $0.5$\% fall in the range of $5\% < \text{GRE} \leq 10\%$. 
The results for $\xi=20$ are comparable. In contrast, the trained CNN with $\xi=10$, yields noticeably poorer results, with $0.2$\% of predictions exceeding a GRE of $10\%$ and one diverged case. 
In the diverged case, negative scaling factors $C(i) < 0$ occurred after $48$ iterations.
This resulted from a gradual amplification of the $\hat{v}_y$ prediction across the Schwarz iterations, 
as small prediction errors in one subdomain propagated as boundary data to neighboring subdomains,
leading to increasingly overestimated recirculation zones in a stenotic segment. 
The resulting negative scaling factors reversed the flow direction, and the deteriorated velocity field was propagated as input data, causing the Schwarz iteration to diverge.

For the purely data-driven CNN models, no diverging predictions occur. 
However, the prediction accuracies are significantly lower than for the flow-rate-conserving USDS.
The best performing model with $\xi=1$ achieves a GRE $\leq 1\%$ for $39.9$\% of the cases, while $4.7$\% exceed a GRE of $5\%$.
The accuracy further deteriorates for the trained CNNs with $\xi=10$ and $\xi=20$. 
For the cases with higher prediction errors, local errors accumulate iteratively,
leading to a gradual overestimation or underestimation of the velocity values.

\Cref{tab:model_comparison_dd_part} provides a detailed comparison of both models for $\xi=1$,
as $\xi=1$ yields the best results for the current problem setting.
The listed metrics are averaged over all predictions falling into each GRE category. {\KL These values are denoted with an overline, for example,  $\overline{\mathrm{GRE}}$ and $\overline{\mathrm{GRE}^*}$. }
First, we note that the component-specific GRE values for the $y$-component are significantly higher than for the $x$-component, as the larger scale $x$-component is currently better learned by the CNN. 
Also, for small flow rates and mild stenoses, $v_y$ approaches zero, leading to higher relative errors.
The flow-rate-conserving USDS requires on average fewer red-black iterations than the purely data-driven USDS. Thus, the global predictions require less computational time.
Note that in the current serial implementation, one red-black iteration  takes approximately $90$\,ms for the $9$ subdomains.
The highest prediction errors for the flow-rate-conserving model occur at very low flow rates ($Q_{\text{norm}} \approx 0.03$), where the velocity magnitudes are small and the CNN lacks sufficient precision.
In these cases, the  $\overline{\mathrm{GRE}^*}$ values are comparable to  $\overline{\mathrm{GRE}}$ , confirming that the error originates from the CNN itself rather than from the Schwarz iteration.
In the remaining categories,  $\overline{\mathrm{GRE}^*}$ is slightly lower than  $\overline{\mathrm{GRE}}$, indicating minor error propagations across the Schwarz iterations.
The data-driven USDS leads to more pronounced differences between  $\overline{\mathrm{GRE}}$  and  $\overline{\mathrm{GRE}^*}$, for example, in the $1$-$5\%$ category,  $\overline{\mathrm{GRE}^*}=0.64\%$ and  $\overline{\mathrm{GRE}}=1.92\%$, 
indicating stronger error accumulation across iterations.
For both models all predictions reach the set tolerance $\epsilon=1.0\cdot 10^{-4}$.

\begin{table}[htbp!]
    \newcommand{\cb}{\cellcolor{blue!10}}
    \newcommand{\frc}{\cellcolor{white!10!orange!40}}    
    \newcommand{\dd}{\cellcolor{white!50!orange!40}}

\centering
\label{tab:global_prediction_xi1}
\resizebox{\textwidth}{!}{%
\begin{tabular}{ll rrrrrrrr}
\toprule
GRE & USDS & \#Count (\%) & \#Converged  & $\overline{\#\mathrm{Red\text{-}black\,it.}}$ & $\overline{Q_{norm}}$ & $\overline{\mathrm{GRE}^*}$ & $\overline{\mathrm{GRE}}$  & $\overline{\mathrm{GRE}}_x$  & $\overline{\mathrm{GRE}}_y$  \\
\midrule
\multirow{2}{*}{$\leq 1\%$}
 &\frc  FR & \frc 644 (83.6\%) & \frc 644  & \frc 10.6 & \frc 0.52 & \frc 0.39\% & \frc 0.53\% & \frc 0.47\% & \frc  9.56\% \\
 & \dd DD & \dd 307 (39.9\%) & \dd 307  & \dd 24.8 & \dd 0.57 & \dd 0.33\% & \dd  0.66\% &  \dd 0.63\% &  \dd  9.62\% \\
\midrule
\multirow{2}{*}{$1$--$5\%$}
 &\frc FR & \frc 122 (15.8\%) & \frc 122  & \frc 12.1 & \frc 0.39 & \frc 1.13\% & \frc 1.51\% & \frc 1.28\% & \frc 20.32\% \\
 & \dd DD & \dd 427 (55.5\%) & \dd 427 & \dd 29.8 & \dd 0.46 & \dd  0.64\%& \dd 1.92\% & \dd 1.87\% & \dd 13.03\% \\
\midrule
\multirow{2}{*}{$5$--$10\%$}
 &\frc FR & \frc   4\ (0.5\%) & \frc   4 & \frc  3.0 & \frc 0.03 & \frc 5.85\% & \frc 5.99\% & \frc 5.11\% & \frc 88.15\% \\
 & \dd DD & \dd  21\ (2.7\%) & \dd  21 & \dd 41.5 & \dd 0.34  & \dd 1.32\% & \dd 6.57\% & \dd 6.48\% & \dd 32.07\% \\
\midrule
\multirow{2}{*}{$10$--$20\%$}
 & FR &  -- & -- &  -- & -- & -- &  -- &  -- &  -- \\
 & \dd DD & \dd  12\ (1.6\%) & \dd  12  & \dd 29.4 & \dd 0.15 & \dd 3.48\% & \dd 13.79\% & \dd 13.65\% & \dd 49.01\% \\
\midrule
\multirow{2}{*}{$>20\%$}
 & FR & -- & -- &  -- & -- &  -- &  --&  -- &  -- \\
 & \dd DD & \dd   3\ (0.4\%) & \dd   3  & \dd 21.3 & \dd 0.04 & \dd  6.26\% & \dd 25.73\% & \dd 25.53\% & \dd 69.69\% \\
\bottomrule
\end{tabular}
}
\caption{ Comparison of global prediction results for the flow-rate-conserving (FR) and data-driven (DD) USDS with $\xi=1$.
For the different GRE categories, the absolute number of predictions (\#Count) belonging to each categoriy is shown, along with the relative percentage (\% of \#Total=$770$),
and the absolute number of converged predictions (\#Converged) is presented. 
Further, the average of the number of red-black iterations ($\overline{\#\mathrm{Red\text{-}black\,it.}}$), $\overline{Q_{\text{norm}}}$, $\overline{\mathrm{GRE}^*}$, achieved $\overline{\mathrm{GRE}}$, $\overline{\mathrm{GRE}}_x$, and $\overline{\mathrm{GRE}}_y$ for predictions falling in each category are listed.}
\label{tab:model_comparison_dd_part}
\end{table}

\begin{figure}[htbp]
    \begin{minipage}[t]{0.49\textwidth}
        \centering
        \includegraphics[width=\textwidth]{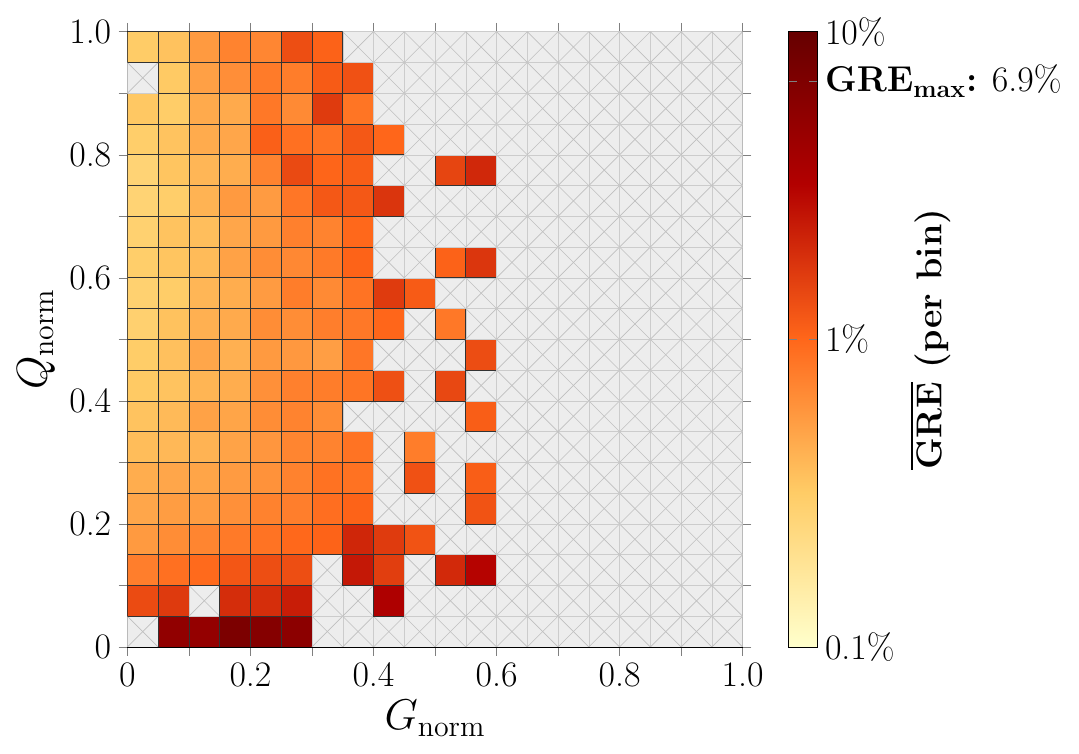}
        \par\smallskip
        \footnotesize\textit{(a) $G_{\text{norm}}$-$Q_{\text{norm}}$ parameter space.}
    \end{minipage}
        \hfill
    \begin{minipage}[t]{0.49\textwidth}
        \centering
        \includegraphics[width=\textwidth]{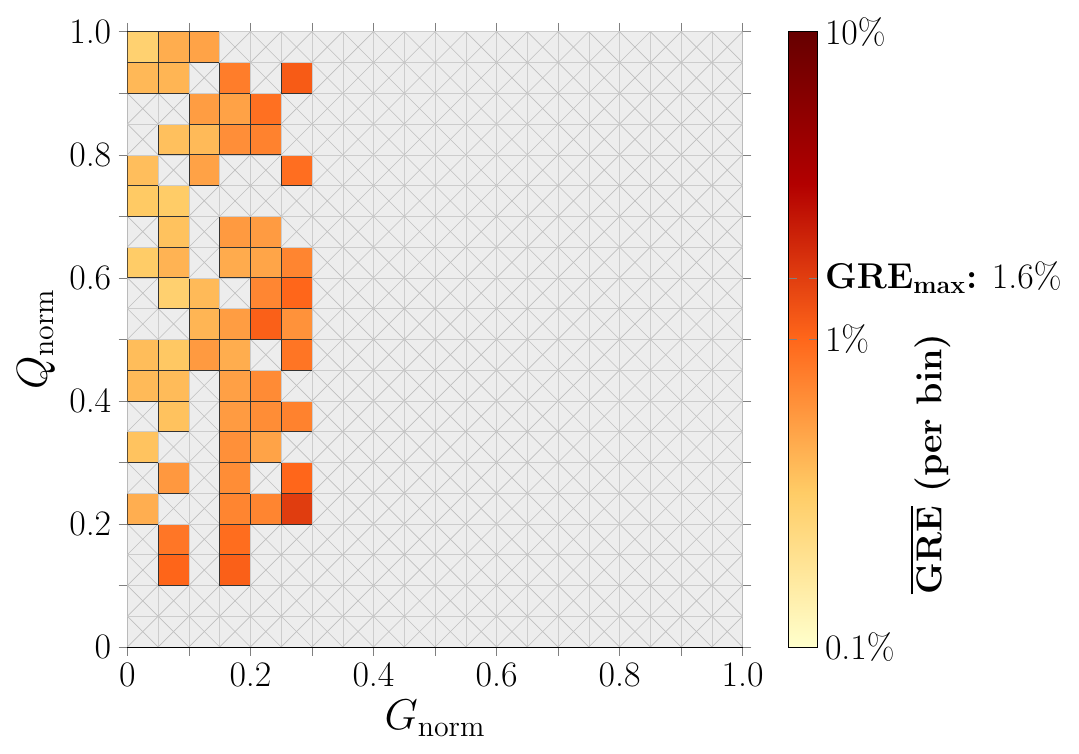}
        \par\smallskip
        \footnotesize\textit{(b) $G_{\text{norm}}$-$Q_{\text{norm}}$ parameter space.}
    \end{minipage}

    \caption{Achieved GRE (averaged per bin) in the $G_{\text{norm}}$-$Q_{\text{norm}}$ parameter space using as USDS the flow-rate-conserving CNN with $\xi=1$. (a) Results for the $770$ problem cases from \cref{sec:dataheatmaps}. (b) Results for $90$ additional test geometries. }
    \label{fig:heatmaps_gre}
\end{figure}

Focusing on the flow-rate-conserving USDS with $\xi=1$, we further analyze the achieved GRE values in the parameter space $G_{\text{norm}}$-$Q_{\text{norm}}$ shown in \cref{fig:heatmaps_gre} (a), where brighter colors indicate a lower average GRE per bin.
We observe that the upper-left region of this plane consistently relates to small errors (GRE $\leq 1\%$), indicating that the method performs well for mild to moderate stenoses with smooth transitions.
Without a theoretical framework, we argue that such systematic analyses are essential to define an admissible parameter range for the usability of the method.
In the $S_{\text{norm}}$-$Q_{\text{norm}}$ space, the results follow a similar trend but do not exhibit as clear a separation between lower and higher error regions as in the $G_{\text{norm}}$-$Q_{\text{norm}}$ space.
This suggests that $G_{\text{norm}}$ is the more informative parameter for determining a feasible application regime.
Based on these results and the data density heatmaps (\cref{fig:heatmaps}), we identify $G_{\text{norm}} \in (0, 0.3]$, $Q_{\text{norm}} \in [0.1, 1.0]$, and $S_{\text{norm}} \in (0, 0.9]$ as parameter ranges for using the trained model as USDS.
To further test this range, we performed $90$ additional predictions with randomly sampled geometries and flow conditions within these bounds.
The results are shown in \cref{fig:heatmaps_gre} (b). We observe that all $90$ predictions converge and achieve a GRE below $2\%$.
Importantly, within the defined range, the flow is viscous-dominated with smooth, moderate stenoses and the relevant flow patterns are well represented in the training data.
Extending the approach beyond these bounds is left for future work; notably, for steeper, more severe stenoses at higher flow velocities, non-unique stationary solutions can occur already at moderate Reynolds numbers and the solution 
becomes sensitive to small geometric perturbations \cite{samuelsson2015}, posing additional challenges for the surrogate model.

\subsection{Transferability to Arterial Domains of Increasing Length}
\label{scalability}
 
Finally, we investigate for the best-performing model from \cref{global_pred_performance}, that is, 
the flow-rate-conserving CNN model with $\xi = 1$, the performance for longer arteries and analyze the convergence behavior in more detail. 
The stenotic region is scaled by the factors $l \in \{2, 4, 8\}$, corresponding to lengths $L_{\text{stenotic}} \in \{2, 4, 8\}$\,[cm] and $14$, $24$, and $44$ subdomains, respectively.
The parameters for the Schwarz algorithm remain the same as in \cref{tab:dd_combined_params}, expcept the stopping tolerance, which is set to $\epsilon=1.0\cdot 10^{-5}$.
In all simulations, the inlet velocity $v^{\max}_{\text{inlet}}=0.3\,[\frac{\text{m}}{\text{s}}]$ is used, such that $Q_{\text{norm}}=0.5$.
Thus, the developed velocity profile for the inlet and outlet conditions is extracted from a single FEM simulation and reused across all tests.

\begin{figure}[htbp]
    \begin{minipage}[t]{0.49\textwidth}
        \centering
        \includegraphics[width=\textwidth]{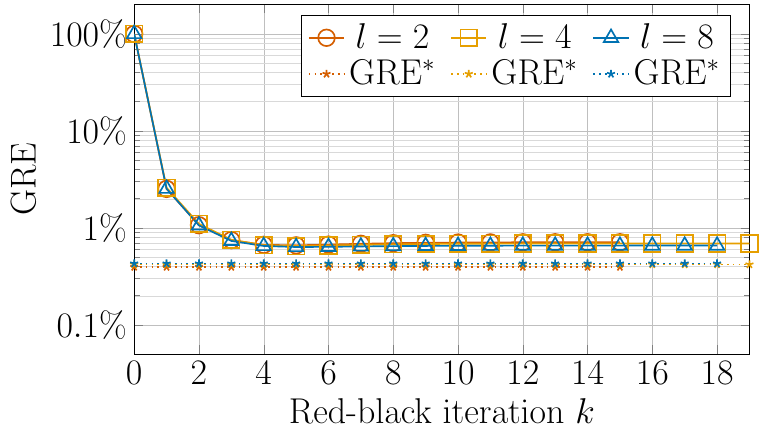}
        \par\smallskip
        \footnotesize\textit{(a) Test Case $1$: Duplicated Stenoses}
    \end{minipage}
    \hfill
    \begin{minipage}[t]{0.49\textwidth}
        \centering
        \includegraphics[width=\textwidth]{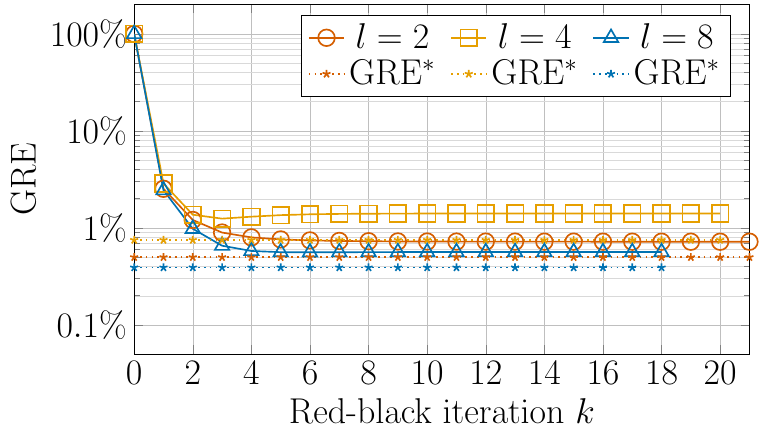}
        \par\smallskip
        \footnotesize\textit{(c) Test Case $2$: Random Stenoses}
    \end{minipage}
        \begin{minipage}[t]{0.49\textwidth}
        \centering
        \includegraphics[width=\textwidth]{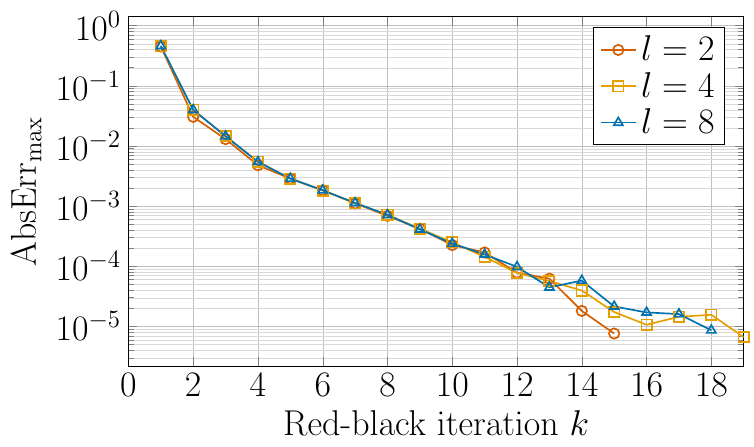}
        \par\smallskip
        \footnotesize\textit{(b) Test Case $1$: Duplicated Stenoses}
    \end{minipage}
    \hfill
    \begin{minipage}[t]{0.49\textwidth}
        \centering
        \includegraphics[width=\textwidth]{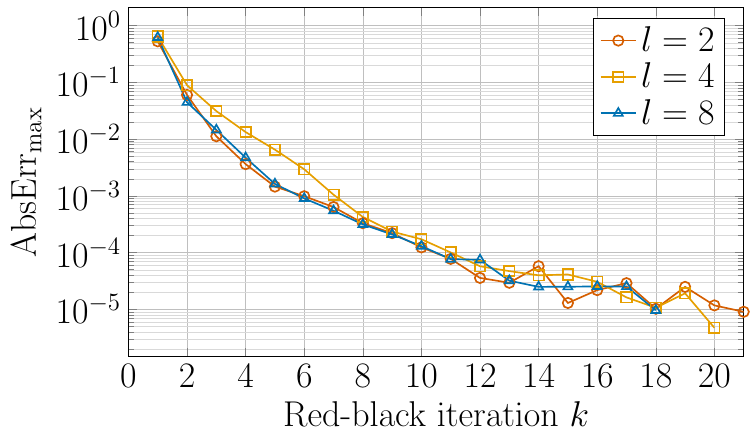}
        \par\smallskip
        \footnotesize\textit{(d) Test Case $2$: Random Stenoses}
    \end{minipage}
    \caption{Convergence behavior of the CNN-Schwarz-Flow method for arteries of increasing length ($l \in \{2, 4, 8\}$) using the flow-rate-conserving CNN with $\xi=1$ as USDS. 
    Top row: GRE over red-black iterations with GRE$^*$ as lower bound. Bottom row: $\text{AbsErr}_{\text{max}}$ over red-black iterations.
     Left column: Duplicated stenotic segments. Right column: Randomly generated stenotic segments.}
    \label{fig:scalability_duplicated_random}
\end{figure}

\paragraph{\textbf{Test Case $1$: Duplicated Stenoses}}
In the first study, the geometries are generated by copying a single stenotic segment ($n_{\text{segment}}=3$, $S_{\text{norm}}=0.49$, $G_{\text{norm}}=0.19$) by the factor $l$,
such that the flow patterns repeat along the arterial length and the dimensionless parameters stay the same.
Note that, as the subdomains overlap, the subdomains in the copied regions do not cover the identical flow patterns due to slight shifts in the subdomain positions.
The results for the three tests in terms of the development of the GRE and stopping criterion, cf. \cref{eq:max:abs_error}, across the iterations  are shown in \cref{fig:scalability_duplicated_random} (a) and (b).
For all three scaling factors, the GRE decreases rapidly during the first five red-black iterations, 
reaching values $<1\%$, and then stagnates. The resulting GRE values are close to the respective GRE$^*$$\approx 0.5\%$ values.
The method converges for all three cases, reaching the tolerance $\epsilon=10^{-5}$ after $15$, $19$, and $18$ iterations for $l=2$, $4$, and $8$, respectively.
The $\text{AbsErr}_{\text{max}}$ decreases approximately linearly on the logarithmic scale during the first iterations but begins to flatten and oscillates slightly below $\epsilon<10^{-4}$ for $l=4$ and $l=8$.
Notably, in all cases, the GRE reaches its final value after approximately eight iterations. 
In the remaining iterations only $\text{AbsErr}_{\text{max}}$ decreases without significantly affecting the GRE.
These observations support the choice of $\epsilon=10^{-4}$ used in the previous section.
As all copied segments share similar flow patterns, this initial study primarily demonstrates that the Schwarz iteration
does not lead to significant error accumulation as the number of subdomains increases.
The method scales in the sense that both the final GRE value and the number of iterations remain comparable for $l = 2$, $4$, and $8$,
even though the number of subdomains increases from $14$ to $44$.
We attribute this, on the one hand, to the parabolic initialization, which provides each interior subdomain with a reasonable starting velocity field, and, on the other hand,
to the flow-rate-conserving constraint layer, which provides global information through the prescribed inlet flow rate to each interior subdomain. 

\begin{figure}[htbp]
    \begin{minipage}[t]{\textwidth}
        \centering
        \includegraphics[width=\textwidth]{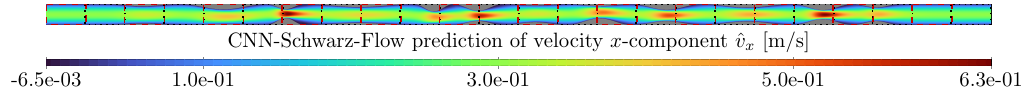}
    \end{minipage}
    \hfill
    \begin{minipage}[t]{\textwidth}
        \centering
        \includegraphics[width=\textwidth]{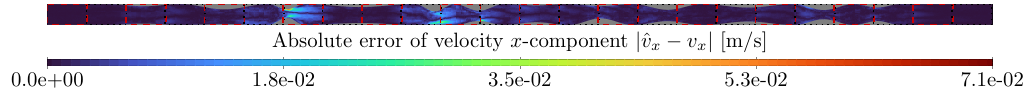}
    \end{minipage}
    \hfill
    \begin{minipage}[t]{\textwidth}
        \centering
        \includegraphics[width=\textwidth]{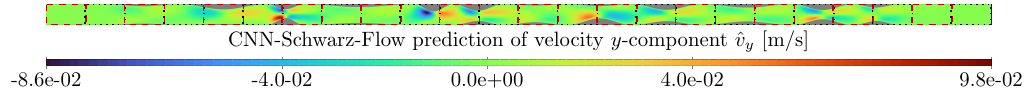}
    \end{minipage}
    \hfill
    \begin{minipage}[t]{\textwidth}
        \centering
        \includegraphics[width=\textwidth]{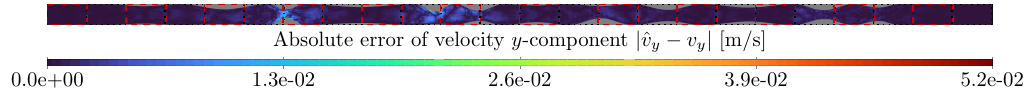}
    \end{minipage}
    \vspace{-19pt}
    \caption{Test geometry for $l=4$ with randomly generated stenotic
segments. The domain is decomposed into $N=24$ subdomains with an overlap of $\delta=2$. The
predicted $x$- and $y$-component using the CNN-Schwarz-Flow approach with the locally trained flow-rate-conserving model for $\xi=1$ after $20$ red-black iterations are visualized. Furthermore, the
absolute errors between the predicted and FEM velocity components are shown.}
    \label{fig:scalability}
\end{figure}

\paragraph{\textbf{Test Case: Random Stenotic Regions}}
In the second study, the geometries consist of randomly generated stenotic segments.
The geometric parameters are as follows; $l=2$: ($n_{\text{segment}}=4$, $S_{\text{norm}} = 0.6$, $G_{\text{norm}}=0.24$), $l=4$: ($n_{\text{segment}}=12$, $S_{\text{norm}}=0.76$, $G_{\text{norm}} = 0.46$), and $l=8$: ($n_{\text{segment}}=17$, $S_{\text{norm}}=0.7$, $G_{\text{norm}}=0.26$).
Among the three cases, the geometry for $l=4$ contains the most severe stenoses and as $G_{\text{norm}}=0.46 > 0.3$ it lies outside our previously defined admissible parameter range; see \cref{fig:heatmaps_gre}.
The results for the three tests are shown in \cref{fig:scalability} (c) and (d).
The GRE$^*$ values now differ between the three tests, reflecting the varying geometric complexity and the corresponding local performance of the USDS for these cases.
The highest GRE$^*\approx 0.7\%$ is observed for $l=4$, that is, for the geometry with the largest $G_{\text{norm}}$.
For $l=2$ and $l=8$, the convergence behavior is comparable to the previous study, with the GRE reaching its plateau within approximately eight iterations and stagnating close to the respective GRE$^*$ values.
In contrast, the GRE for $l=4$ decreases during the first three iterations, but then increases again before eventually stagnating at a GRE$\approx1.5\%$.
Compared to the other tests, the $\text{AbsErr}_{\text{max}}$ shows larger changes between consecutive iterations during the initial eight iterations in this case.
An analysis of the error in each subdomain reveals that the subdomain with the most severe stenosis ($G_{\text{norm}}=0.46$) exhibits the largest local errors; see \cref{fig:scalability}.
While errors in other subdomains are decreasing, the error in this specific subdomain starts to increase after four iterations, until it stagnates after ten iterations.
Concluding, the results demonstrate that the method is transferable to longer arteries with varying stenotic regions.
The method converges more reliably and achieves smaller errors in the test cases, where the local flow configurations remain within the empirically determined admissible parameter range.

\section{Conclusion}
\label{sec:conclusion}

The proposed CNN-Schwarz-Flow algorithm enables the prediction of stationary flow fields in domains of arbitrary sizes, in our case,
arteries of arbitrary length, by using locally trained CNN surrogate models within an alternating overlapping Schwarz method.
The integration of a physics-aware constraint layer enforcing flow rate conservation significantly improves the global prediction accuracy and reduces iteration counts
compared to a purely data-driven approach.
We emphasize that defining admissible parameter bounds for the local geometric and flow configurations covered by the training data is essential, as the CNN is an inexact subdomain solver.
Within these bounds, the method converges to global predictions with small errors and scales to longer arteries with increasing numbers of subdomains.
While the CNN inference time per subdomain is on the order of milliseconds, the computational cost of generating training data, training the CNN, and validating the model is significant.
An advantage over classical solvers is therefore expected only if the trained model is reused across many problem instances of larger dimension, for example, in patient-specific simulations with varying inflow conditions.
This is not specific for our model, but applies principally to all learned surrogate or reduced-order models.

\section*{Acknowledgments}
\paragraph{Computational Ressources} 
This research was conducted using the Fujitsu \\PRIMEHPC FX1000 (Wisteria/BDEC-01 (Odyssey)) and the Fujitsu PRIMERGY\\ GX2570 equipped with NVIDIA A100 Tensorcore (Wisteria/BDEC-01 (Aquarius)) at the Information Technology Center, the University of Tokyo under support by ``Joint Usage/Research Center for Interdisciplinary Large-scale Information Infrastructures (JHPCN) (jh240029, jh250027)''.
The authors gratefully acknowledge the scientific support and HPC resources provided by the Erlangen National High Performance Computing Center (NHR@FAU) of the Friedrich-Alexander-Universität Erlangen-Nürnberg (FAU) under the NHR project k105be. 
NHR funding is provided by federal and Bavarian state authorities. NHR@FAU hardware is partially funded by the German Research Foundation (DFG) – 440719683.
\paragraph{Data} We gratefully acknowledge the group of G. Gompper at the Institute for Advanced Simulation (IAS-2) at Forschungszentrum Jülich, including D. Fedosov
 and A. Topuz, 
 for providing us the viscosity data (\cref{tab:viscosity_params}).
\paragraph{Code Developement} The initial code framework for subdomain and global prediction of
 stationary flows around circular objects was implemented by S. Hatayama 
for his master's thesis "Steady flow prediction across multiple regions using deep learning and boundary exchange" (Original text in Japanese; translated into English), supervised by T. Shimokawabe at the Department of Electrical Engineering and Information Systems, University of Tokyo.
 We gratefully acknowledge S. Hatayama for his contributions.
\paragraph{Declaration of Generative AI and AI-assisted technologies in the
writing and coding process}
While preparing this work, N. Kubicki used AI tools including DeepL, Claude (Sonnet 3.7), 
and GitHub Copilot (Tool=GPT-4o/4.1) to assist with grammar and spelling corrections, improve writing style and advance code development.
\bibliographystyle{siamplain}

\bibliography{references/references}

@book{kundu2024fluid,
  author = {Pijush K. Kundu and Ira M. Cohen},
  title = {Fluid Mechanics},
  edition = {4th},
  publisher = {Elsevier},
  year = {2010},
  address = {New York}
  }

@article{samuelsson2015,
  author  = {Samuelsson, J. and Tammisola, O. and Juniper, M. P.},
  title   = {Breaking axi-symmetry in stenotic flow lowers the critical transition {R}eynolds number},
  journal = {Phys. Fluids},
  volume  = {27},
  number  = {10},
  eid   = {104103},
  year    = {2015},
  doi     = {10.1063/1.4934530},
}

@article{RAISSI2019686,
title = {Physics-informed neural networks: A deep learning framework for solving forward and inverse problems involving nonlinear partial differential equations},
journal = {J. Comput. Phys.},
volume = {378},
pages = {686--707},
year = {2019},
doi = {10.1016/j.jcp.2018.10.045},
author = {M. Raissi and P. Perdikaris and G.E. Karniadakis}
}

@techreport{USDpt:SciML,
title        = {Brochure on {B}asic {R}esearch {N}eeds for {S}cientific {M}achine {L}earning: {C}ore {T}echnologies for {A}rtificial {I}ntelligence},
author = {Baker, Nathan and Alexander, Frank and Bremer, Timo and Hagberg, Aric and Kevrekidis, Yannis and Najm, Habib and Parashar, Manish and Patra, Abani and Sethian, James and Wild, Stefan and Willcox, Karen},
url = {https://www.osti.gov/biblio/1484362}, 
institution  = {USDOE Office of Science (SC)},
place = {US},
year = {2018}
}

@techreport{frey2001medit,
  author = {Pascal Frey},
  title = {MEDIT: An interactive Mesh visualization Software},
  institution = {INRIA},
  number = {RT-0253},
  year = {2001}
  }

@article{McGreivy2024_weakbaseline,
  author = {McGreivy, Nick and Hakim, Ammar},
  title = {Weak baselines and reporting biases lead to overoptimism in machine learning for fluid-related partial differential equations},
  journal = {Nat. Mach. Intell.},
  volume = {6},
  pages = {1256--1269},
  year = {2024},
  doi = {10.1038/s42256-024-00897-5},
}

@misc{review_ml_cfd_wang,
author = {Wang, Haixin and Cao, Yadi and Huang, Zijie and Liu, Yuxuan and Hu, Peiyan and Luo, Xiao and Song, Zezheng and Zhao, Wanjia and Liu, Jilin and Sun, Jinan and Zhang, Shikun and Wei, Long and Wang, Yue and Wu, Tailin and Ma, Zhi-Ming and Sun, Yizhou},
title = {Recent Advances on Machine Learning for Computational Fluid Dynamics: A Survey},
note = {preprint, arXiv:2408.12171, 2024}
}

@book{CardiovascularMathematics,
  title     = "Cardiovascular Mathematics: Modeling and simulation of the circulatory system",
  author    = "Formaggia, Luca  and Quarteroni, Alfio and Veneziani,  Alessandro ",
  year      = "2009",
  publisher = "Springer",
  address   = "Milan",
  doi="10.1007/978-88-470-1152-6",
}

@book{irgens_rheology,
  title     = "Rheology and Non-Newtonian Fluids",
  author    = "Fridtjov Irgens",
  year      = "2014",
  publisher = "Springer",
  address   = "Cham, Switzerland",
  doi="10.1007/978-3-319-01053-3",
}

@article{Bessonov2016_MethodsBloodFlow,
author = {Bessonov, N. and Sequeira, A. and Simakov, S. and Vassilevskii, Yu. and Volpert, V.},
doi = {10.1051/mmnp/201611101},
journal = {		Math. Model. Nat. Phenom.},
pages = {1--25},
title = {Methods of Blood Flow Modelling},
volume = {11},
year = {2016}
}

@incollection{Fasano2017,
author="Fasano, Antonio
and Sequeira, Ad{\'e}lia",
title="Hemorheology and Hemodynamics",
bookTitle="Hemomath: The Mathematics of Blood ",
year="2017",
publisher="Springer, Cham",
pages="1--77",
address="Switzerland",
}

@incollection{Sequeira2018,
author="Sequeira, A.",
editor="Farina, Angiolo
and Mikeli{\'{c}}, Andro
and Rosso, Fabio",
title="Hemorheology: Non-{N}ewtonian Constitutive Models for Blood Flow Simulations",
bookTitle="Non-Newtonian Fluid Mechanics and Complex Flows",
year="2018",
publisher="Springer, Cham",
address="Switzerland",
pages="1--44"
}

@Misc{Yadan2019Hydra,
  author =       {Omry Yadan},
  title =        {Hydra - A framework for elegantly configuring complex applications},
  year =         {2019},
  note         = {GitHub repository},
  url =          {https://github.com/facebookresearch/hydra}
}

@inproceedings{pytorch,
author = {Paszke, Adam and Gross, Sam and Massa, Francisco and Lerer, Adam and Bradbury, James and Chanan, Gregory and Killeen, Trevor and Lin, Zeming and Gimelshein, Natalia and Antiga, Luca and Desmaison, Alban and K\"{o}pf, Andreas and Yang, Edward and DeVito, Zach and Raison, Martin and Tejani, Alykhan and Chilamkurthy, Sasank and Steiner, Benoit and Fang, Lu and Bai, Junjie and Chintala, Soumith},
title = {PyTorch: an imperative style, high-performance deep learning library},
year = {2019},
publisher = {Curran Associates Inc.},
address = {New York},
booktitle = {Proceedings of the 33rd International Conference on Neural Information Processing Systems},
articleno = {721},
}

@misc{schlomer_pygmsh_2022,
  author       = {Nico Schlömer},
  title        = {Pygmsh: A Python frontend for Gmsh},
  year         = {2022},
  doi          = {10.5281/zenodo.1173105},
  url          = {https://github.com/nschloe/pygmsh},
  note         = {GitHub repository}
}

@book{vtkBook,
  author    = "Will Schroeder and Ken Martin and Bill Lorensen",
  title     = "The Visualization Toolkit: An Object-Oriented Approach to 3D Graphics",
  edition = {4th},
  publisher = "Kitware",
  year      = "2006",
  address   = "New York",
}

@inproceedings{kingma2017adammethodstochasticoptimization,
      title={Adam: A Method for Stochastic Optimization}, 
      author={Diederik P. Kingma and Jimmy Ba},
      year={2015},
      booktitle = {International Conference for Learning Representations},
      address = {San Diego},
      note = {arXiv:1412.6980 [cs.LG]}
}

@misc{horovod,
  Author = {Alexander Sergeev and Mike Del Balso},
  Title = {Horovod: fast and easy distributed deep learning in {TensorFlow}},
  note      = {preprint, 	arXiv:1802.05799 [cs.LG], 2018}

}

@misc{feddlib,
  author       = { {FEDDLib Team} },
  title        = {{FEDDL}ib ({F}inite Element and Domain Decomposition Library)},
  year         = {2025},
  note         = {GitHub repository},
  url          = {https://github.com/FEDDLib/FEDDLib}
  }

@misc{MATLAB,
author = {{The MathWorks Inc.}},
year = {2022},
title = {MATLAB Version: 9.13.0 (R2022b)},
publisher = {The MathWorks Inc.},
address = {Natick, Massachusetts, United States},
url = {https://www.mathworks.com}
}

@misc{plotnn,
  author    = {Haris Iqbal},
  title     = {PlotNeuralNet},
  year      = {2020},
  howpublished = {\url{https://github.com/HarisIqbal88/PlotNeuralNet}},
  note      = {GitHub repository}
}

@incollection{lecun:1989:CNN,
  title={Generalization and network design strategies},
  author={LeCun, Yann},
  booktitle = {Connectionism in Perspective},
  editor = { R. Pfeifer and Z. Schreter and F. Fogelman and L. Steels},
  year={1989},
  publisher={Elsevier},
  address = {Zurich}
}

@book{Goodfellow:2016:DL,
	Author = {Goodfellow, Ian and Bengio, Yoshua and Courville, Aaron},
	Publisher = {MIT Press Cambridge},
	Title = {Deep learning},
	year = {2016},
  address = {Cambridge, Massachusetts},
  note={\url{http://www.deeplearningbook.org}},
}

@inproceedings{guo_cnn,
author = {Guo, Xiaoxiao and Li, Wei and Iorio, Francesco},
year = {2016},
pages = {481--490},
title = {Convolutional Neural Networks for Steady Flow Approximation},
publisher = {Association for Computing Machinery},
address = {New York},
booktitle = {Proceedings of the 22nd ACM SIGKDD International Conference on Knowledge Discovery and Data Mining}
}

@article{eichinger_cnn,
author = {Eichinger, Matthias and Heinlein, Alexander and Klawonn, Axel},
year = {2022},
pages = {235--255},
title = {Surrogate convolutional neural network models for steady computational fluid dynamics simulations},
volume = {56},
journal = {Electron. Trans. Numer. Anal.},
doi = {10.1553/etna_vol56s235}
}

@incollection{prechelt2002early,
  title={Early stopping-but when?},
  author={Prechelt, Lutz},
  booktitle={Neural Networks: Tricks of the trade},
  pages={55--69},
  year={2002},
  publisher={Springer},
  address = {Cham, Switzerland},
  doi = {10.1007/3-540-49430-8_3},
  volume    = {1524},
  series = {Lecture Notes Comput. Sci.},
}

@article{article_rana_cnn_dd,
author = {Rana, Pratip and Weigand, Timothy and Pilkiewicz, Kevin and Mayo, Michael},
year = {2024},
eid = {23080},
title = {A scalable convolutional neural network approach to fluid flow prediction in complex environments},
volume = {14},
journal = {Sci. Rep.},
doi = {10.1038/s41598-024-73529-y}
}

@article{Kim2019DeepFluids,
  title={Deep Fluids: A Generative Network for Parameterized Fluid Simulations},
  author={Byungsoo Kim and Vinicius C. Azevedo and Nils Thuerey and Theodore Kim and Markus Gross and Barbara Solenthaler},
  journal={Comput. Graph. Forum},
  volume={38},
  number={2},
  pages={59--70},
  year={2019},
  doi={10.1111/cgf.13619}
}

@article{PhysRevFluids.8.014604,
  title = {Embedding hard physical constraints in neural network coarse-graining of three-dimensional turbulence},
  author = {Mohan, Arvind T. and Lubbers, Nicholas and Chertkov, Misha and Livescu, Daniel},
  journal = {Phys. Rev. Fluids},
  volume = {8},
  eid = {014604},
  year = {2023},
  doi = {10.1103/PhysRevFluids.8.014604},
}

@article{weather_cnn_downsampling,
author = {Harder, Paula and Hernandez-Garcia, Alex and Ramesh, Venkatesh and Yang, Qidong and Sattegeri, Prasanna and Szwarcman, Daniela and Watson, Campbell D. and Rolnick, David},
title = {Hard-constrained deep learning for climate downscaling},
year = {2023},
volume = {24},
pages = {1--40},
journal = {J. Mach. Learn. Res.},
}

@article{article_beucler_constraint,
author = {Beucler, Tom and Pritchard, Michael and Rasp, Stephan and Ott, Jordan and Baldi, Pierre and Gentine, Pierre},
title = {Enforcing Analytic Constraints in Neural Networks Emulating Physical Systems},
year = {2021},
eid = {098302},
volume = {126},
journal = {Phys. Rev. Lett.},
doi = {10.1103/PhysRevLett.126.098302}
}

@article{rheometer,
author = {Fedosov, Dmitry and Karniadakis, George and Caswell, Bruce},
year = {2010},
eid= {144103},
title = {Steady shear rheometry of dissipative particle dynamics models of polymer fluids in reverse Poiseuille flow},
volume = {132},
journal = {J. Chem. Phys.},
doi = {10.1063/1.3366658}
}

@article{article_he,
author = {He, Xin and Neytcheva, Maya and Vuik, C.},
year = {2015},
pages = {33-58},
title = {On Preconditioning of Incompressible Non-{N}ewtonian Flow Problems},
volume = {33},
journal = {J. Comput. Math.},
doi = {10.4208/jcm.1407-m4486}
}

@book{volker,
title={Finite Element Methods for Incompressible Flow Problems},
author={Volker John},
year={2016},
publisher={Springer},
address = {Berlin},
series    ={Springer Series in Comput. Math.},
volume    = {51},
doi ={10.1007/978-3-319-45750-5}
}

@article{pachecognf,
title = "On outflow boundary conditions in finite element simulations of non-{N}ewtonian internal flows",
author = "Pacheco, {Douglas Ramalho Queiroz} and M{\"u}ller, {Thomas Stephan} and Olaf Steinbach and G{\"u}nter Brenn",
year = "2021",
doi = "10.51375/IJCVSE.2021.1.6",
volume = "1",
journal = "Int. J. Comput. Vis. Sci. Engrg.",
eid = "6"
}

@article{new_paper_schwarz_pde_neural_operators,
title = "A Learning-based Domain Decomposition Method",
author = "Wu, Rui  and Kovachki, Nikola  and  Liu, Burigede",
year = "2026",
doi = "10.1016/j.cma.2026.118799",
volume = "453",
journal = "Comput. Methods Appl. Mech. Engrg.",
eid = "118799"
}

@misc{preprint,
author = {Klaes, Simon and Klawonn, Axel and Kubicki, Natalie and Lanser, Martin and Nakajima, Kengo and Shimokawabe, Takashi and Weber, Janine},
title = {A Flow-rate-conserving CNN-based Domain Decomposition Method for Blood Flow Simulations},
note = {preprint,	arXiv:2509.15900 [math.NA], 2025}
}

@article{article_dd_ml_survey,
author = {Klawonn, Axel and Lanser, Martin and Weber, Janine},
year = {2024},
title = {Machine learning and domain decomposition methods - a survey},
volume = {1},
journal = {Comput. Sci. Engrg.},
doi = {10.1007/s44207-024-00003-y},
eid={2}
}

@article{Wang_2022,
   title={Mosaic flows: A transferable deep learning framework for solving {PDE}s on unseen domains},
   volume={389},
   DOI={10.1016/j.cma.2021.114424},
   journal={Comput. Methods Appl. Mech. Engrg.},
   author={Wang, Hengjie and Planas, Robert and Chandramowlishwaran, Aparna and Bostanabad, Ramin},
   year={2022},
   eid={114424} 
   }

@inproceedings{feeney2023breakingboundariesdistributeddomain,
author = {Feeney, Arthur and Li, Zitong and Bostanabad, Ramin and Chandramowlishwaran, Aparna},
title = {Breaking Boundaries: Distributed Domain Decomposition with Scalable Physics-Informed Neural {PDE} Solvers},
year = {2023},
publisher = {Association for Computing Machinery},
booktitle = {Proceedings of the International Conference for High Performance Computing, Networking, Storage and Analysis (SC '23)},
eid = {80},
numpages = {1--15},
location = {Denver, USA}}

@book{book_domaindecomposition,
  author    = {Toselli, Andrea and Widlund, Olof},
  title     = {Domain Decomposition Methods -- Algorithms and Theory},
  series    ={Springer Series in Comput. Math.},
  volume    = {34},
  year      = {2005},
  doi       = {10.1007/b137868},
  publisher = {Springer},
  address   = {Berlin}
}

@book{dd_quarteroni,
    author    = {Quarteroni, Alfio and Valli, Alberto},
    title     = {Domain Decomposition Methods for Partial Differential Equations},
    publisher = {Oxford Uni. Press},
    year      = {1999},
    doi       = {10.1093/oso/9780198501787.001.0001},
    address   = {Oxford}
}

@misc{monolithic_preconditioner_lea,
author = {Heinlein, Alexander and Klawonn, Axel and Knepper, Jascha and Saßmannshausen, Lea},
title = {{M}onolithic and Block Overlapping {S}chwarz Preconditioners for the Incompressible {N}avier--{S}tokes Equations},
note      = {preprint, 	arXiv:2506.16179 [math.NA], 2025}
}

@phdthesis{hochmuth,
            year = {2020},
           title = {Parallel Overlapping Schwarz Preconditioners for Incompressible Fluid Flow and Fluid-Structure Interaction Problems},
          school = {University of Cologne},
          author = {Christian Hochmuth},
             url = {https://kups.ub.uni-koeln.de/11345/},
            type = {{Ph.D.} thesis},
        address  = {Cologne, Germany}
}

@incollection{Frosch,
author    = {Heinlein, Alexander and Klawonn, Axel and Rajamanickam,  Sivasankaran and Rheinbach, Oliver},
doi       = {10.1007/978-3-030-56750-7_19},
booktitle = {Domain Decomposition Methods in Science and Engineering XXV},
series    = {Lect. Notes Comput. Sci. Eng.},
volume    = {138},
publisher = {Springer},
address   = {Cham, Switzerland},
year      = {2020},
pages     = {176--184},
title = {{FROS}ch: A Fast And Robust Overlapping {S}chwarz Domain Decomposition Preconditioner Based on {X}petra in {T}rilinos}
}

@misc{mayr2025trilinosenablingscientificcomputing,
      title={{T}rilinos: {E}nabling Scientific Computing Across Diverse Hardware Architectures at Scale}, 
      author={Matthias Mayr and Alexander Heinlein and Christian Glusa and Siva Rajamanickam and Maarten Arnst and Roscoe Bartlett and Luc Berger-Vergiat and Erik Boman and Karen Devine and Graham Harper and Michael Heroux and Mark Hoemmen and Jonathan Hu and Brian Kelley and Kyungjoo Kim and Drew P. Kouri and Paul Kuberry and Kim Liegeois and Curtis C. Ober and Roger Pawlowski and Carl Pearson and Mauro Perego and Eric Phipps and Denis Ridzal and Nathan V. Roberts and Christopher Siefert and Heidi Thornquist and Romin Tomasetti and Christian R. Trott and Raymond S. Tuminaro and James M. Willenbring and Michael M. Wolf and Ichitaro Yamazaki},
      note={preprint, arXiv:2503.08126 [cs.MS], 2025}
}

\end{document}